\input lanlmac
\def\href#1#2{{#2}}
\def\hhref#1{{#1}}
\input epsf.tex

\overfullrule=0mm

\newcount\figno
\figno=0
\def\fig#1#2#3{
\par\begingroup\parindent=0pt\leftskip=1cm\rightskip=1cm\parindent=0pt
\baselineskip=11pt
\global\advance\figno by 1
\midinsert
\epsfxsize=#3
\centerline{\epsfbox{#2}}
\vskip 12pt
{\bf Fig.\ \the\figno:} #1\par
\endinsert\endgroup\par
}
\def\figlabel#1{\xdef#1{\the\figno}}
\def\encadremath#1{\vbox{\hrule\hbox{\vrule\kern8pt\vbox{\kern8pt
\hbox{$\displaystyle #1$}\kern8pt}
\kern8pt\vrule}\hrule}}


\def\IR{\relax{\rm I\kern-.18em R}}
\font\cmss=cmss10 \font\cmsss=cmss10 at 7pt

\font\cmss=cmss10 \font\cmsss=cmss10 at 7pt
\def\IZ{\relax\ifmmode\mathchoice
{\hbox{\cmss Z\kern-.4em Z}}{\hbox{\cmss Z\kern-.4em Z}}
{\lower.9pt\hbox{\cmsss Z\kern-.4em Z}}
{\lower1.2pt\hbox{\cmsss Z\kern-.4em Z}}\else{\cmss Z\kern-.4em Z}\fi}
\def\IN{\relax{\rm I\kern-.18em N}}
\def\b{\circ}
\def\n{\bullet}

\def\gbbbb{\Gamma_4^{\hbox{$\scriptstyle \b \b$}\kern -8.2pt
\raise -4pt \hbox{$\scriptstyle \b \b$}}}
\def\gnnnn{\Gamma_4^{\hbox{$\scriptstyle \n \n$}\kern -8.2pt  
\raise -4pt \hbox{$\scriptstyle \n \n$}}}
\def\gnnnnnn{\Gamma_6^{\hbox{$\scriptstyle \n \n \n$}\kern -12.3pt
\raise -4pt \hbox{$\scriptstyle \n \n \n$}}}
\def\gbbbbbb{\Gamma_6^{\hbox{$\scriptstyle \b \b \b$}\kern -12.3pt
\raise -4pt \hbox{$\scriptstyle \b \b \b$}}}
\def\gbbbbc{\Gamma_{4\, c}^{\hbox{$\scriptstyle \b \b$}\kern -8.2pt
\raise -4pt \hbox{$\scriptstyle \b \b$}}}
\def\gnnnnc{\Gamma_{4\, c}^{\hbox{$\scriptstyle \n \n$}\kern -8.2pt
\raise -4pt \hbox{$\scriptstyle \n \n$}}}
\def\Rbud#1{{\cal R}_{#1}^{-\kern-1.5pt\blacktriangleright}}
\def\Rleaf#1{{\cal R}_{#1}^{-\kern-1.5pt\vartriangleright}}
\def\Rbudb#1{{\cal R}_{#1}^{\circ\kern-1.5pt-\kern-1.5pt\blacktriangleright}}
\def\Rleafb#1{{\cal R}_{#1}^{\circ\kern-1.5pt-\kern-1.5pt\vartriangleright}}
\def\Rbudn#1{{\cal R}_{#1}^{\bullet\kern-1.5pt-\kern-1.5pt\blacktriangleright}}
\def\Rleafn#1{{\cal R}_{#1}^{\bullet\kern-1.5pt-\kern-1.5pt\vartriangleright}}
\def\Wleaf#1{{\cal W}_{#1}^{-\kern-1.5pt\vartriangleright}}
\def\Cleaf{{\cal C}^{-\kern-1.5pt\vartriangleright}}
\def\Cbud{{\cal C}^{-\kern-1.5pt\blacktriangleright}}
\def\Crleaf{{\cal C}^{-\kern-1.5pt\circledR}}


\Title{\vbox{\hsize=3.truecm \hbox{SPhT/04-060}}}
{{\vbox {
\bigskip
\centerline{Planar maps as labeled mobiles}
}}}
\bigskip
\centerline{J. Bouttier\foot{bouttier@spht.saclay.cea.fr}, 
P. Di Francesco\foot{philippe@spht.saclay.cea.fr} and
E. Guitter\foot{guitter@spht.saclay.cea.fr}}
\medskip
\centerline{ \it Service de Physique Th\'eorique, CEA/DSM/SPhT}
\centerline{ \it Unit\'e de recherche associ\'ee au CNRS}
\centerline{ \it CEA/Saclay}
\centerline{ \it 91191 Gif sur Yvette Cedex, France}
\bigskip
\noindent 
We extend Schaeffer's bijection between rooted quadrangulations
and well-labeled trees to the general case of Eulerian planar maps 
with prescribed face valences to obtain a bijection with a new class 
of labeled trees, which we call mobiles. 
Our bijection covers all the classes of maps previously
enumerated by either the two-matrix model used by physicists 
or by the bijection with blossom trees used by combinatorists.
Our bijection reduces the enumeration of maps to that, much
simpler, of mobiles and moreover keeps track of the geodesic distance 
within the initial maps via the mobiles' labels.
Generating functions for mobiles are shown to obey systems of
algebraic recursion relations.

\Date{05/04}

\nref\TUT{W. Tutte,
{\it A Census of Planar Maps}, Canad. J. of Math. 
{\bf 15} (1963) 249-271.}
\nref\BIPZ{E. Br\'ezin, C. Itzykson, G. Parisi and J.-B. Zuber, {\it Planar
Diagrams}, Comm. Math. Phys. {\bf 59} (1978) 35-51.}
\nref\CORV{R. Cori and B. Vauquelin, {\it Planar maps are well labeled trees},
Canad. J. Math. {\bf 33 (5)} (1981) 1023-1042.}
\nref\ARQUES{D. Arqu\`es, {\it Les hypercartes planaires sont des arbres 
tr\`es bien \'etiquet\'es}, Discr. Math. {\bf 58}(1) (1986) 11-24.}  
\nref\SCHth{G. Schaeffer, {\it Conjugaison d'arbres
et cartes combinatoires al\'eatoires}, PhD Thesis, Universit\'e 
Bordeaux I (1998); see 
\hhref{http://www.loria.fr/\~{}schaeffe/index-en.html}.}
\nref\GEOD{J. Bouttier, P. Di Francesco and E. Guitter, {\it Geodesic
distance in planar graphs}, Nucl. Phys. {\bf B663}[FS] (2003) 535-567, 
arXiv:cond-mat/0303272.}
\nref\LALLER{J. Bouttier, P. Di Francesco and E. Guitter, {\it Random
trees between two walls: Exact partition function}, J. Phys. A: Math. Gen.
{\bf 36} (2003) 12349-12366, arXiv:cond-mat/0306602.} 
\nref\ONEWALL{J. Bouttier, P. Di Francesco and E. Guitter, {\it Statistics
of planar maps viewed from a vertex: a study via labeled trees},
Nucl. Phys. {\bf B675}[FS] (2003) 631-660, arXiv:cond-mat/0307606.}
\nref\CS{P. Chassaing and G. Schaeffer, {\it Random Planar Lattices and 
Integrated SuperBrownian Excursion}, 
Probability Theory and Related Fields {\bf 128(2)} (2004) 161-212, 
arXiv:math.CO/0205226.}
\nref\CHASDUR{P. Chassaing and B. Durhuus, {\it Statistical Hausdorff 
dimension of labelled trees and quadrangulations}, 
arXiv:math.PR/0311532.}
\nref\MARMO{J. F. Marckert and A. Mokkadem, {\it Limit of normalized 
quadrangulations: the Brownian map}, arXiv:math.PR/0403398.} 
\nref\BECAN{E. Bender and E. Canfield, {\it The number of degree-restricted 
rooted maps on the sphere}, SIAM J. Discrete Math. {\bf 7}(1) (1994) 9-15.}
\nref\CENSUS{J. Bouttier, P. Di Francesco and E. Guitter, {\it Census of planar
maps: from the one-matrix model solution to a combinatorial proof},
Nucl. Phys. {\bf B645}[PM] (2002) 477-499, arXiv:cond-mat/0207682.}
\nref\SCH{G. Schaeffer, {\it Bijective census and random
generation of Eulerian planar maps}, Electronic
Journal of Combinatorics, vol. {\bf 4} (1997) R20.}
\nref\TUTSLI{W. Tutte, 
{\it A Census of slicings},
Canad. J. of Math. {\bf 14} (1962) 708-722.}
\nref\BOUKA{D. Boulatov and V. Kazakov, {\it The Ising model
on a random planar lattice: the structure of the phase 
transition and the exact critical exponents}, Phys. Lett. {\bf B186} (1987) 379-384.}
\nref\DOU{M. Douglas, {\it The two-matrix model}, in 
{\it Random Surfaces and Quantum Gravity}, O. Alvarez,
E. Marinari and P. Windey eds., NATO ASI Series {B:} Physics Vol. {\bf 262} (1990).} 
\nref\BMS{M. Bousquet-M\'elou and G. Schaeffer,{\it The degree distribution
in bipartite planar maps: application to the Ising model},
arXiv:math.CO/0211070.}
\nref\CONST{M. Bousquet-M\'elou and G. Schaeffer,
{\it Enumeration of planar constellations}, Adv. in Applied Math.,
{\bf 24} (2000) 337-368.}
\nref\HUR{A. Hurwitz, {\it \"Uber Riemann'sche Fl\"achen mit gegebenen 
Verzweignungspunkten}, Math. Ann. {\bf 39} (1891) 1-60.}

\newsec{Preliminaries}

\subsec{Introduction}

Maps, like graphs and trees, are fundamental combinatorial objects that
appear in many different areas of mathematics, computer science, and
theoretical or statistical physics. The groundwork for the enumerative theory
of planar maps was laid in the 60's by Tutte \TUT, who enumerated
maps of some particular classes, with some remarkably simple results. 
For example, the number of rooted planar maps with a given number $n$ 
of edges is $2 (2n)! 3^n \over n!  (n+2)!$.
For many other classes of maps, no similar closed general formula is known
but the corresponding generating functions can
be shown to obey algebraic equations. Such results are obtained
using, for instance, the formalism developed by Tutte (the so-called
recursive decomposition) or the powerful method of matrix integrals
developed by physicists \BIPZ.

However the proofs involved there are rather of an indirect,
non-constructive nature. Of greater mathematical beauty are bijective
proofs, where the purpose is to establish one-to-one correspondences
between classes of maps and some simpler sets whose enumeration is
obvious. Such correspondences were first found in the 80's by Cori and
Vauquelin \CORV\ and later Arqu\`es \ARQUES, but the real development of
this subject came with the thesis of Schaeffer \SCHth, who was able to
rederive most of Tutte's results through new bijective algorithms cutting
maps into trees. It was then realized that these bijections gave insight
into some detailed information on the intrinsic geometry of maps such as
the geodesic distance, with interesting applications both in statistical
physics [\xref\GEOD-\xref\ONEWALL] and probability theory
[\xref\CS-\xref\MARMO].

\subsec{Known results, aim of the paper}

We now introduce a few definitions in order to recall some known
results. A {\it planar map} is an embedding of a connected graph 
(which may have loops and/or parallel edges) in the
sphere such that the edges are non-intersecting open curves and connected only
at their extremities (vertices). The complement of the graph is
a disjoint union of simply connected domains (faces). In
enumerative theory, two maps differing by an homeomorphism (bicontinuous
bijection) of the sphere are considered equivalent, hence the number of
maps with a finite number of edges is also finite. Most known results deal
with {\it rooted} maps, that is maps having a distinguished oriented edge 
(the {\it root}).

In a map, the {\it valence} or {\it degree} of a face or vertex is
the number of its incident edges\foot{These are counted
with multiplicity : more precisely we count the number of incident edge
sides (resp.\ half-edges).}. Maps with only 4-valent faces are
called {\it quadrangulations} which are dual to tetravalent maps
where all vertices have degree 4. This class of maps is the one for which
bijections are the simplest. There are actually two of them \SCHth: 
one between tetravalent maps and so-called blossom trees, 
the other between quadrangulations and well-labeled trees.

So far, generalizations to wider classes of maps were obtained mainly in
terms of blossom trees.  A first extension consists of enumerating
maps with prescribed vertex valences. The corresponding generating
functions are easily derived by considering the general ``one-matrix
model", or by recursive decomposition \BECAN, while the corresponding
bijective proof involves blossom trees with subtle charge
constraints \CENSUS. This is drastically simplified in the case of maps
with only vertices of even valence, corresponding to even potentials in
the matrix model formulation. The corresponding blossom trees \SCH\ then
have a simple characterization which makes it possible to re-derive Tutte's 
compact formulas for the numbers of maps with prescribed (even) vertex valences
\TUTSLI.

More generally, one may as well try to enumerate the {\it bipartite}, i.e.
vertex-bicolored, maps with prescribed vertex valences of either color,
which corresponds to the general ``two-matrix model". This problem has
many physical applications, including, for instance, the celebrated Ising
model on random lattices \BOUKA\ as well as a whole range of
multicritical theories corresponding to minimal models of CFT coupled to
2D quantum gravity \DOU.  The bijective enumeration via blossom trees was
found by Bousquet-M\'elou and Schaeffer \BMS. This indeed extends the
previous case, as arbitrary maps are equivalent to bipartite maps with
only two-valent black vertices. Another interesting subcase is that of
$p$-constellations with only $p$-valent black vertices and white vertices
with valences that are multiples of $p$ \CONST. This corresponds to the 
most general
situation where explicit compact formulas \HUR\ are known for the numbers
of maps with prescribed vertex valences. To complete the picture, note
that $2$-constellations are equivalent to maps with even vertex
valences.

On the other (dual) front, well-labeled trees appeared so far only in
correspondence with quadrangulations \SCHth\ \CS\ and Eulerian
(face-bicolored) triangulations \ONEWALL\foot{Quadrangulations and
Eulerian triangulations are respectively related to the original bijection
of Cori and Vauquelin \CORV\ and to that of Arqu\`es \ARQUES, which both
rely on an encoding of maps in terms of permutations.}. Here the labels on
the trees correspond to a labeling of the vertices of the original map
by their geodesic distance from a fixed origin vertex, which gives
access to a host of results on the geometry of maps.  This is to be
contrasted with the blossom tree approach, where only two-point
correlations (i.e. generating functions for maps with two points at a
fixed geodesic distance) are within reach.

The aim of this paper is to extend this construction to the general case
of Eulerian planar maps with prescribed face valences (dual to the bipartite
planar maps with prescribed vertex valences considered in Ref.\ \BMS) 
by introducing generalized
labeled trees, which we call mobiles. In terms of generating
functions, we reproduce and generalize the algebraic recursion relations
of \GEOD\ and provide a new combinatorial interpretation for them.

\subsec{Plan of the paper}

The paper is organized as follows. For the sake of clarity, we begin in
Sect.\ 2 by the simpler case of planar maps with even face valences. In
Sect.\ 2.1, we associate to each such map with a distinguised vertex 
a mobile whose
properties are further characterized. In Sect.\ 2.2, we display the
inverse construction, proved in detail in Sect.\ 2.3. Finally, in
Sect.\ 2.4 we derive equations for the corresponding generating functions, 
which determine a single function $R_n$ for maps with two points at a 
geodesic distance less than $n$.

In Sect.\ 3, we turn to the general case of Eulerian maps with prescribed
face valences : we construct the corresponding mobiles in Sect.\ 3.1, and
we exhibit the inverse construction in Sect.\ 3.2. Generating functions 
are studied in Sect.\ 3.3.

Sect.\ 4 is devoted to the illustration of the sub-cases of
$p$-constellations and arbitrary maps, where mobiles may be simplified. A
few concluding remarks on interesting applications are made in
Sect.\ 5.

\newsec{The case of planar maps with prescribed even face degrees}

\subsec{From maps to mobiles}

We start from a planar map $\cal M$ whose faces all have even degrees.
Equivalently, this amounts to requiring that the map be {\it bipartite},
namely that its vertices may be partitioned into two sets
so that no two adjacent vertices belong to the same set.
On this map, we distinguish a vertex as {\it the origin} and label all 
the vertices
of the map by their {\it geodesic distance} to this origin,
i.e. the length of any shortest path from the origin to that
vertex. With this definition, the origin is the only vertex labeled
$0$, its nearest neighbors are all labeled $1$, ..., and all the
labels are non-negative integers. Moreover, any two adjacent vertices
have labels differing by at most one, and because of the bipartite nature of 
the map, these labels cannot be equal.
\fig{A typical configuration of labels around a face of $\cal M$
(delimited by dashed edges), 
here of valence $2k=10$. Adjacent labels differ by $\pm 1$. 
We add a black unlabeled vertex at the center of the face and connect it 
(via solid edges) to the $k=5$ labeled vertices immediately followed clockwise 
by a smaller label.}{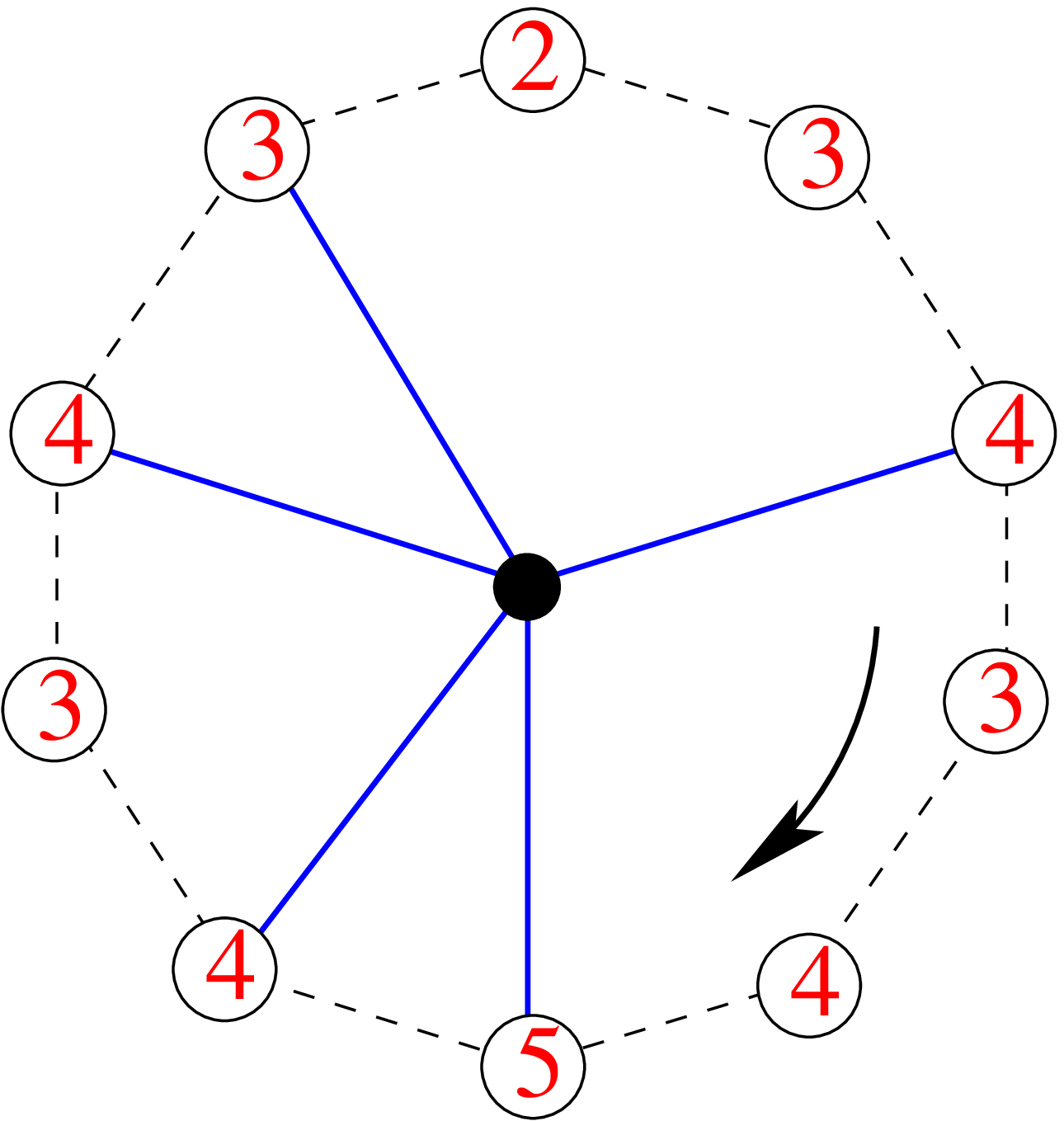}{4.cm}
\figlabel\select

Our construction takes place independently within each face of $\cal M$.
Given a face of degree $2k$, the adjacent vertex labels
read in clockwise direction around the face form a cyclic sequence
with increments of $\pm 1$. Among these $2k$ vertices, we select 
the $k$ ones immediately followed by vertices with smaller labels
(see figure \select\ for an example).
We add a new (unlabeled) vertex at the center of the face and
connect it within the face by $k$ new non-intersecting edges to 
the $k$ selected labeled vertices.
After completing this construction within each face, we remove all the 
edges of the original map. We also erase the origin vertex as
by construction it becomes isolated. We are left with a map
$\cal T$ with two types of vertices: unlabeled ones in correspondence 
with the faces of $\cal M$, and labeled ones consisting of all
the vertices of $\cal M$ except the origin. The edges of $\cal T$
connect only vertices of different types.
\fig{Proof by contradiction that $\cal T$ has no cycle.
Assuming the existence of such a cycle, we pick a labeled vertex with minimal
label $n$ on the cycle. The neighborhoods of its two adjacent unlabeled
vertices on the cycle show the existence of a vertex labeled $n-1$ on
each side of the cycle. A path on $\cal M$ from the origin to one of 
these vertices must cross the cycle at a labeled vertex with label
at most $n-1$, a contradiction.}{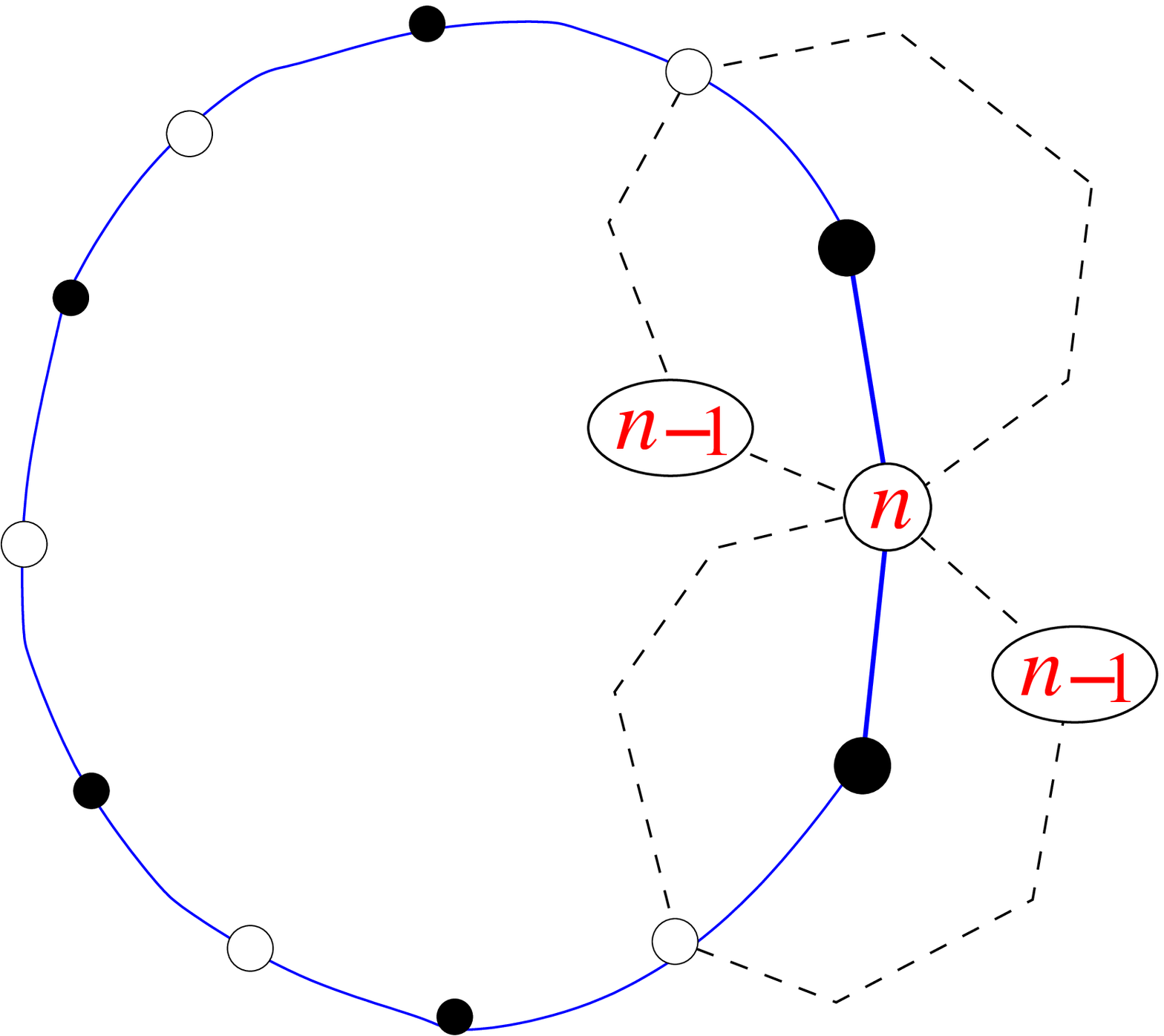}{5.cm}
\figlabel\surround
We now show that $\cal T$ is a plane tree. We first show that 
it contains no cycle and then that it has a single connected component.
The first statement is proved by contradiction. Assume that $\cal T$
contains a cycle, namely a closed path separating the plane into two regions.
We call the interior of the cycle the region not containing the former
origin. We pick a labeled vertex on the cycle whose label,
say $n$, is minimal along the cycle. By examining the neighborhood of
its two adjacent unlabeled vertices on the cycle (see figure \surround), 
we conclude that there is a vertex labeled $n-1$ in the interior
of the cycle. This is a contradiction: a geodesic path on $\cal M$
from the origin to this vertex must intersect the cycle at a labeled
vertex with label $j\le n-1< n$.

\fig{A typical example of a planar map $\cal M$ (a) with one 2-valent,
three 4-valent and two 6-valent faces. Having selected an origin (labeled $0$),
all the vertices are naturally labeled by their geodesic distance to that
origin. We perform (b) the construction of figure \select\ in each face
of $\cal M$ (with, as usual, the opposite conventions for the external 
face due to the representation in the plane). Erasing the original edges
of $\cal M$ as well as its origin, we are left (c) with a tree 
$\cal T$.}{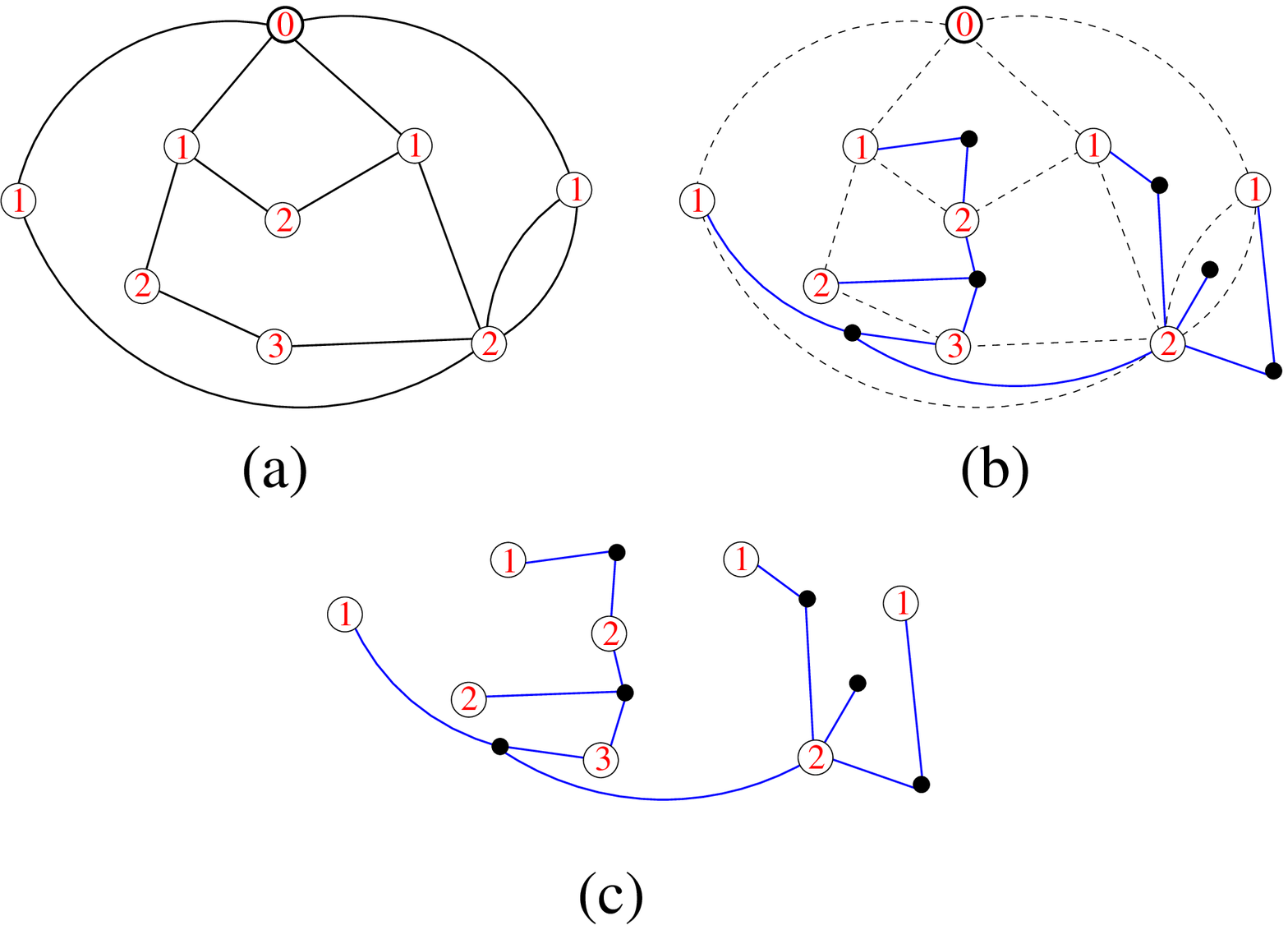}{12.cm}
\figlabel\maptotree
Having shown that $\cal T$ has no cycle, it is a forest made of 
$c$ trees. Let $V$, $F$ and $E$ denote respectively the numbers
of vertices, faces and edges of $\cal M$, obeying the Euler relation
$V-E+F=2$. The total number of vertices of $\cal T$ is $F+V-1$.
The number of edges in $\cal T$ reads $\sum_k k F_{2k}$ where 
$F_{2k}$ is the number of $2k$-valent faces of $\cal M$. This
is nothing but $E$ as clearly $\sum_k 2k F_{2k}=2E$.
By noting that each of the $c$ trees has one more vertex than edge,
we deduce that $c=(F+V-1)-E=1$. $\cal T$ is therefore a tree.
An illustration of the complete construction of $\cal T$ starting
from $\cal M$ is shown in figure \maptotree.

The labels of ${\cal T}$ have the following property: 
\item{(P)} for each unlabeled vertex $v$, the labels $n$ and $m$ of
two labeled vertices adjacent to $v$ and consecutive in clockwise
direction satisfy $m\geq n-1$.

We define a {\it mobile} as a plane tree obeying the following rules:
\item{(i)} its vertices are of two types, unlabeled ones and 
labeled ones carrying integer labels,
\item{(ii)} each edge connects a labeled to an unlabeled vertex,
\item{(iii)} the labels obey the property (P) above.
\par
\noindent If the mobile moreover satisfies the additional rule:
\item{(iv)} all labels are strictly positive and there is at least 
one vertex with label $1$,
\par
\noindent then it will be said {\it well-labeled}.
The reader may check that figure \maptotree-(c) represents 
a well-labeled mobile.

The above construction maps each bipartite planar map with 
an origin into a well-labeled mobile.
This mapping turns out to be a bijection whose inverse is exhibited 
in the next section.

\subsec{Converse construction}

\fig{A typical example of a well-labeled mobile $\cal T$ (a). For each 
unlabeled vertex of $\cal T$, the successive labels of its adjacent 
labeled vertices cannot decrease by more than $1$ clockwise. We add an 
extra origin vertex labeled $0$ and connect (b) each labeled corner to 
its successor (dashed arrows). Erasing all unlabeled vertices of $\cal T$
and their adjacent edges, we are left (c) with a bipartite planar map
$\cal M$. Moreover, the labels simply encode the geodesic distance
from the origin to the vertices.}{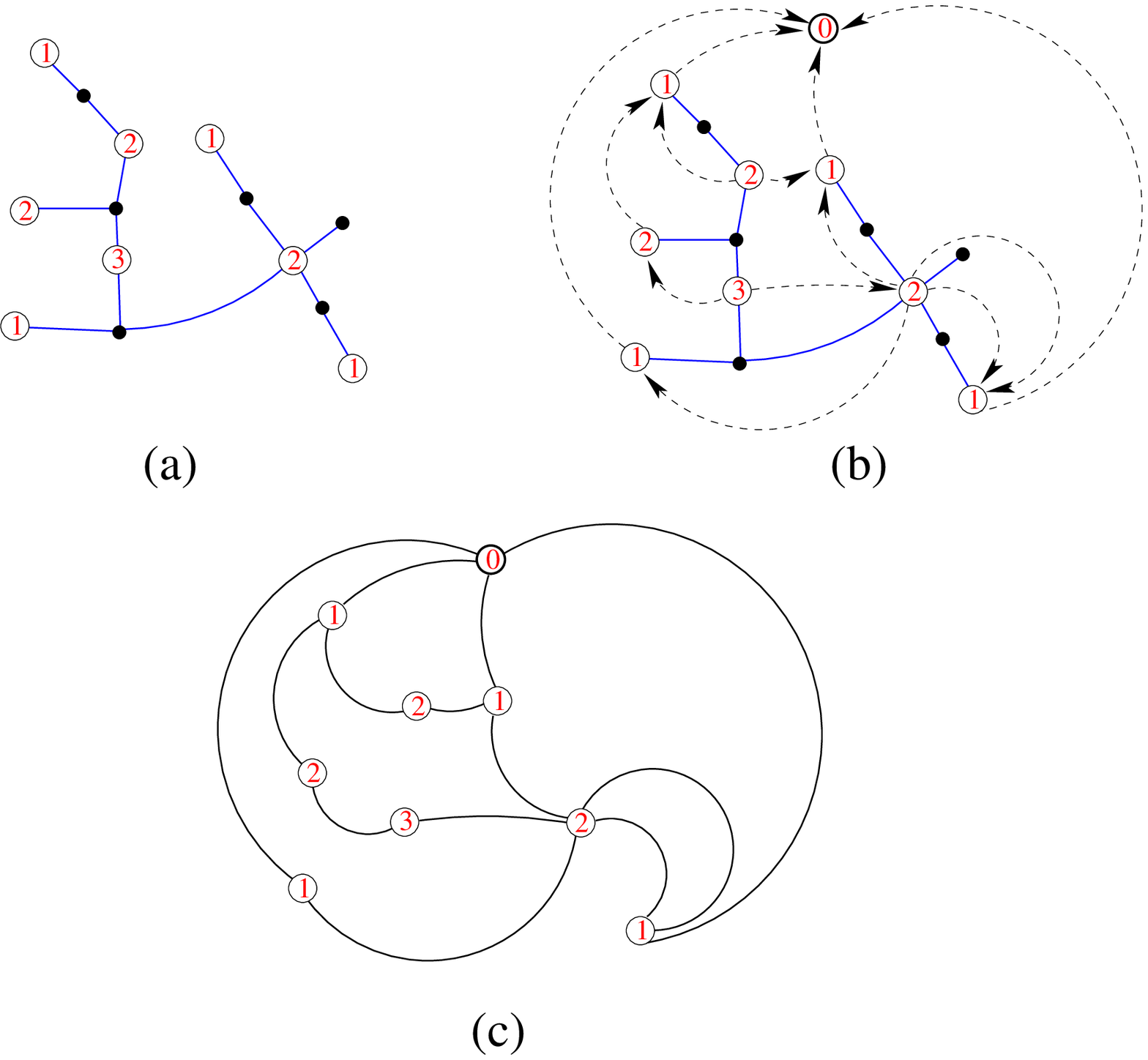}{12.cm}
\figlabel\treetomap
We start from a well-labeled mobile $\cal T$. 
A {\it corner} of $\cal T$ is 
a sector with apex at a labeled vertex of $\cal T$ and delimited by two
consecutive edges around this vertex. We label each corner by the
label of its apex. To each corner $C$ with label $n\geq 2$, 
we associate its successor $s(C)$ defined as the first encountered 
corner with label $n-1$ when going clockwise around the tree 
(see figure \treetomap-(b)). 
The existence of a successor is ensured by the property (P). Indeed,
at each step in the sequence of corners read clockwise around $\cal T$, 
the label may decrease by at most $1$; hence between the corner $C$ and
a corner with label $1$, all labels between $1$ and $n$ must be
present. 

We construct the map $\cal M$ associated with $\cal T$ by first
drawing an edge between each corner with label $n\geq 2$ and
its successor within the external face of $\cal T$ and in such
a way that no two edges intersect. This can be done due
to the nested structure around $\cal T$ of corners and their successors, 
namely that if a corner $C'$ lies strictly between a corner $C$ and 
its successor $s(C)$, then $s(C')$ lies between $C'$ and $s(C)$ (with
possibly $s(C')=s(C)$). Again this is a consequence of the 
property (P). We next add an origin vertex labeled $0$ in the
external face and view the unique sector around this isolated point 
as the successor of all corners labeled $1$, which we therefore also 
connect to the origin via non-crossing edges. This is possible because each 
corner has its successor before or at the first encountered corner 
labeled $1$; hence all corners labeled $1$ are incident to the external 
face.  Finally we erase all unlabeled vertices and their adjacent edges. 
The result is a map $\cal M$ with an origin, which is connected  
because each vertex is connected to the origin via a chain of successors. 
It is moreover bipartite because the parity of labels alternate between 
adjacent 
vertices. We may forget about labels because these are nothing but the geodesic
distances from the vertices of $\cal M$ to the origin. 
Indeed, since the labels on vertices adjacent
in $\cal M$ differ by exactly $1$, the geodesic distance from a vertex
to the origin is larger than or equal to its label $n$, and a chain 
of successors provides a geodesic of length $n$. Figure \treetomap\ displays
an example of construction of the map $\cal M$ starting from a
well-labeled mobile $\cal T$. 

In the next section, we argue that the above construction is indeed
the inverse of that presented in section 2.1.

\subsec{Proof of the bijection}

In order to prove that the two previous constructions are inverse
of one another, we have to show successively the two following assertions:
\item{(1)} starting with a bipartite planar map $\cal M$ with an origin
and constructing its associated mobile $\cal T$ as in section 2.1,
the construction of section 2.2 carried out on this mobile $\cal T$
retrieves $\cal M$,
\item{(2)} starting with a well-labeled mobile $\cal T$
and constructing its associated map $\cal M$ as in section 2.2,
the construction of section 2.1 carried out on this map $\cal M$
retrieves $\cal T$.
\par

\fig{Proof that any edge $e$ of $\cal M$ is restored in the construction
of section 2.2. Apart from the trivial case where the origin is
adjacent to $e$, the union of $e$ and $\cal T$ has two faces. We
follow (a) counterclockwise the contour path $\cal C$ bordering the face not 
containing the origin and pick a corner with smallest label $m$
(represented here with an outgoing arrow). If
this corner were not the last one on $\cal C$ (adjacent to $e$), 
it would be followed by an unlabeled vertex $v$ and by the selection rules
of figure 1 around $v$, we would
deduce the existence of another corner with label $m-1$ inside the
face. As the face does not contain the origin, the corner $m-1$ is different
from the origin hence it lies on the contour $\cal C$, a contradiction. 
This implies (b) that the end corner of $\cal C$ is the successor of its
starting corner, hence that $e$ is restored.}{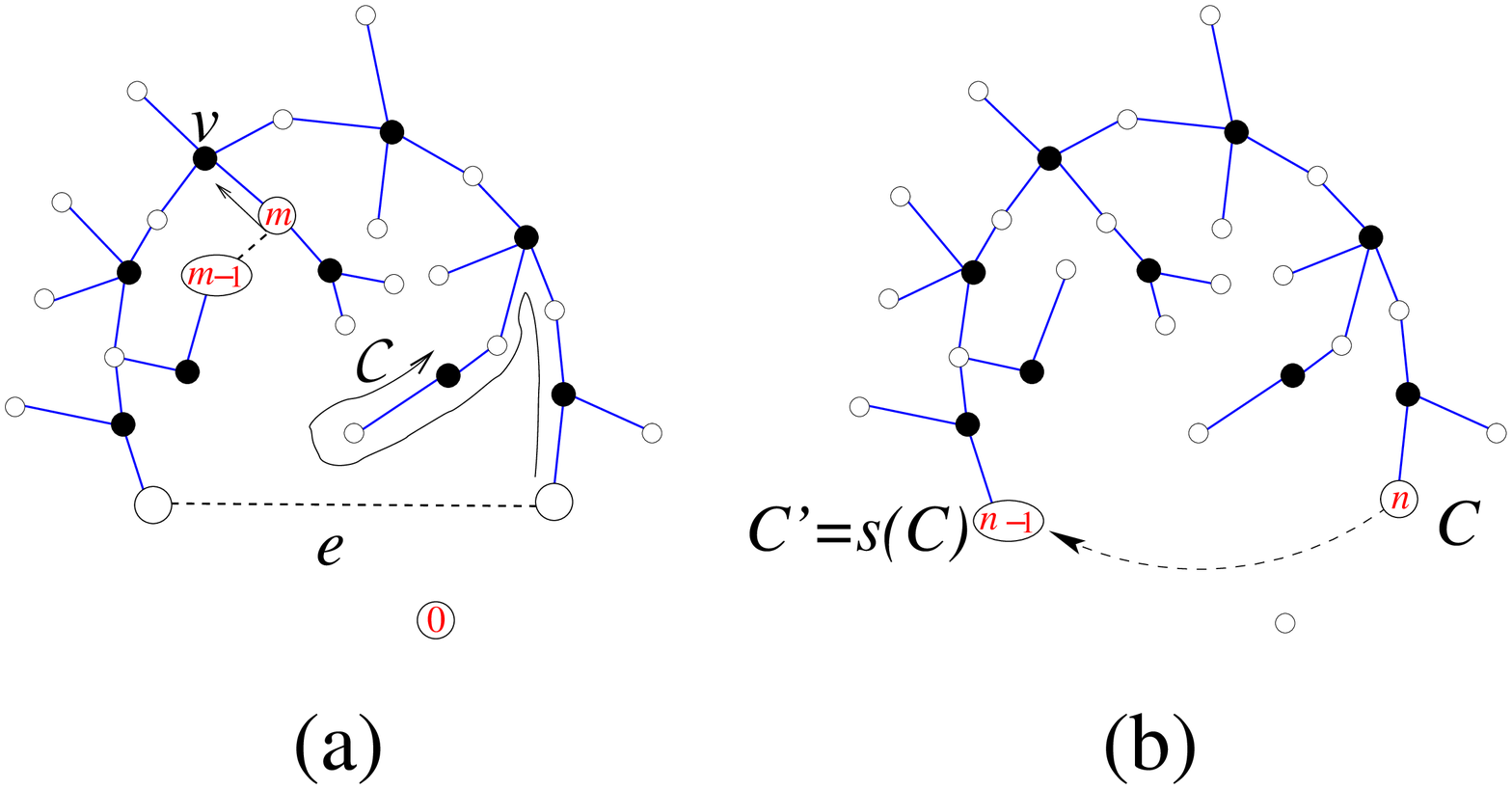}{13.cm}
\figlabel\bijecone

To prove (1), we consider a bipartite map $\cal M$ with an origin and
its associated mobile $\cal T$. We only have to prove that the
construction of section 2.2 restores exactly the edges of $\cal M$, as the
vertices of $\cal M$ are the labeled vertices in $\cal T$ or the origin.
We consider an edge $e$ of $\cal M$ connecting two vertices labeled
$n$ and $n-1$, say. In $\cal T$, the extremities of $e$ point to two
corners of the union ${\cal T}\cup \{0\}$ of the tree $\cal T$ and the
origin vertex $\{0\}$. We denote by $C$ the corner with label $n$
and by $C'$ that with label $n-1$: we have to show that $C'$ is the successor
of $C$. If $n=1$, this is automatic by construction. If $n\geq 2$, both
$C$ and $C'$ are true corners of $\cal T$. We consider the graph obtained
by adding to $\cal T$ the edge $e$: it has exactly $2$ faces, one of which
does not contain the origin. We follow the contour path $\cal C$ going
counterclockwise around this face from one extremity of $e$ to the other.
Let $m$ denote the smallest corner label on $\cal C$. Assume by
contradiction that this label is attained at a corner which is not the
last one on $\cal C$: as shown in figure \bijecone, we deduce the
existence of a corner on $\cal C$ with label $m-1<m$. Hence $m$ is only
attained at the last corner on $\cal C$, which is therefore $C'$ with
label $m=n-1$, while the first corner is $C$: as all labels in between are
strictly larger than $n-1$, this proves that $C'$ is the successor of $C$,
and thus $e$ is restored in the construction of section 2.2. A simple
counting argument shows that the number of corners in ${\cal T}$ is exactly the
number of edges in ${\cal M}$; hence the construction of section 2.2 does
not create edges other than those of ${\cal M}$.

\fig{The generic configuration (a) of successors of two consecutive vertices
adjacent to a given unlabeled vertex of $\cal T$. Each face of $\cal M$ is
obtained (b) as the union around an unlabeled vertex of all such 
configurations. As the labels decrease by one along the arrows, we
immediately see that the labels of $\cal T$ are re-selected by the
construction of section 2.1 acting on $\cal M$.}{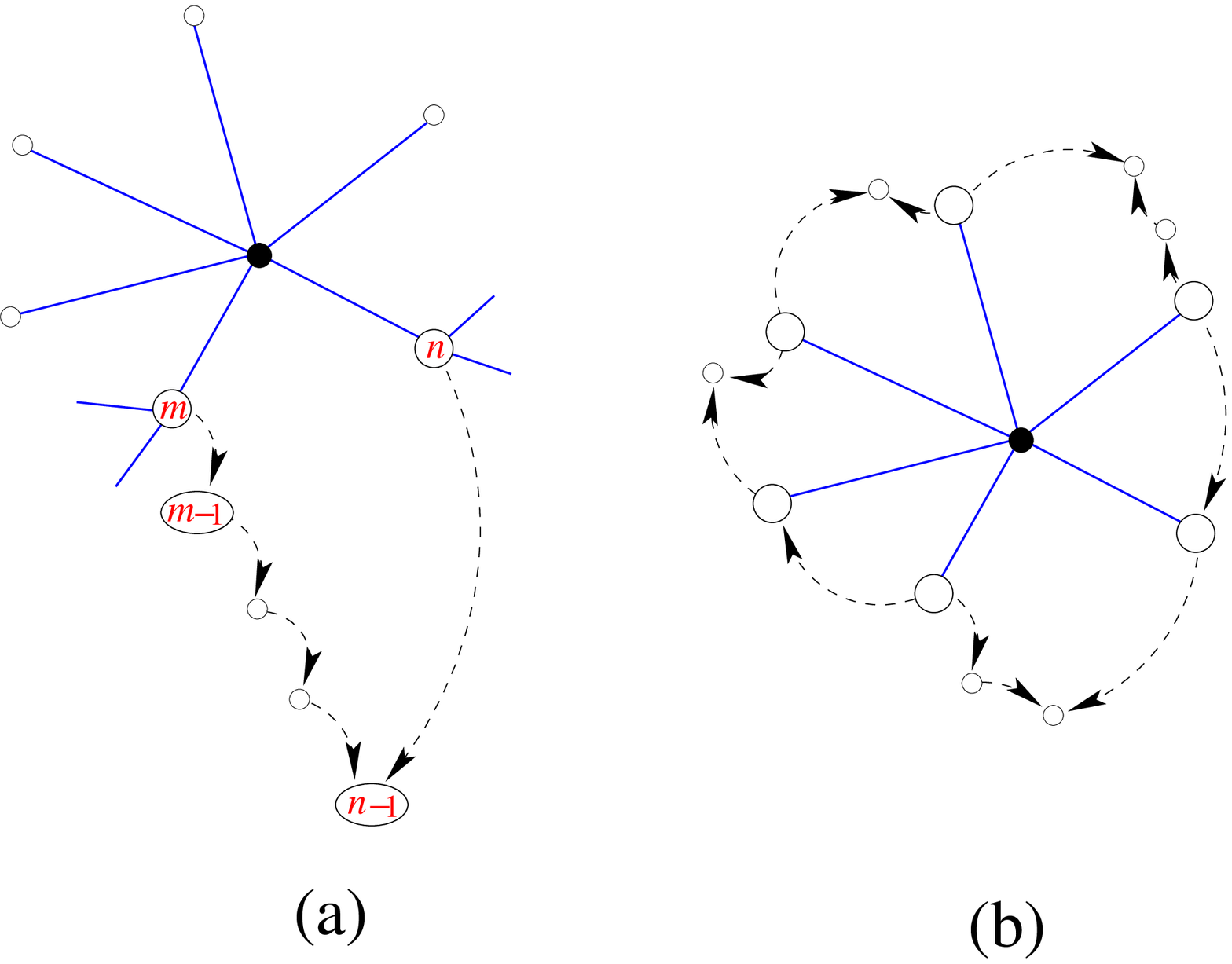}{10.cm}
\figlabel\bijectwo
To prove (2), we consider a well-labeled mobile $\cal T$
and its associated map $\cal M$. We now characterize the faces of 
$\cal M$. We first show that each face of $\cal M$ contains exactly
one unlabeled vertex of $\cal T$, and then that the rules of
section 2.1 inside this face select precisely the edges of $\cal T$
incident to that vertex. Start with an unlabeled vertex of $\cal T$
and consider two clockwise consecutive adjacent corners as in 
figure \bijectwo-(a), with labels $n$ and $m$ (with $m\geq n-1$ by definition 
of a mobile). Then, the successor of the corner labeled $n$  
belongs to the sequence of successors of the corner labeled $m$.
This delimits a region of the plane which contains no other
vertex or edge of $\cal M$. The union of these regions for all 
the corners adjacent to the unlabeled vertex at hand forms a face 
of $\cal M$ containing no other unlabeled vertex. 
All the faces of $\cal M$ are obtained this way, as shown by
a counting argument: the numbers of edges in $\cal T$ and $\cal M$
are equal (and both equal to the number of labeled corners in $\cal T$),
and there is one more vertex in $\cal M$ as labeled vertices in $\cal T$.
Using Euler's relation both for $\cal M$ and $\cal T$  
shows that there are as many faces in $\cal M$ than unlabeled
vertices in $\cal T$. Finally, as apparent from figure \bijectwo-(b), 
the selection rules
of section 2.1 select precisely the vertices originally connected
in $\cal T$ to the same unlabeled vertex. 

\subsec{Generating functions}

The constructions of Sects.\ 2.1 and 2.2 establish a bijection between, on 
the one hand, bipartite planar maps with an origin vertex, and, on the other 
hand, well-labeled mobiles. Enumeration is simpler for rooted mobiles, which 
enjoy recursive properties. More precisely, a rooted mobile has a distinguished 
corner, which in terms of maps corresponds to distinguishing an edge. 
We may moreover attach weights $g_{2k}$ per $k$-valent unlabeled vertex of 
the mobile, which amounts to a weight $g_{2k}$ per $2k$-valent face of the 
planar maps.
\fig{Top: the decomposition of a rooted mobile with root label $n$ 
(generating function $R_n$) into a number of mobiles with univalent root
vertex with the same label (generating function $L_n$) leads to
eq. (2.1). Bottom: A mobile with univalent root vertex with label $n$
is in turn decomposed according to the clockwise neighborhood of
the adjacent unlabeled vertex, in correspondence with walks (see
oval). Each descending step $i\to i-1$ on the walk gives rise
to a factor $R_i$ for the associated rooted sub-mobile, leading
to eq. (2.2).}{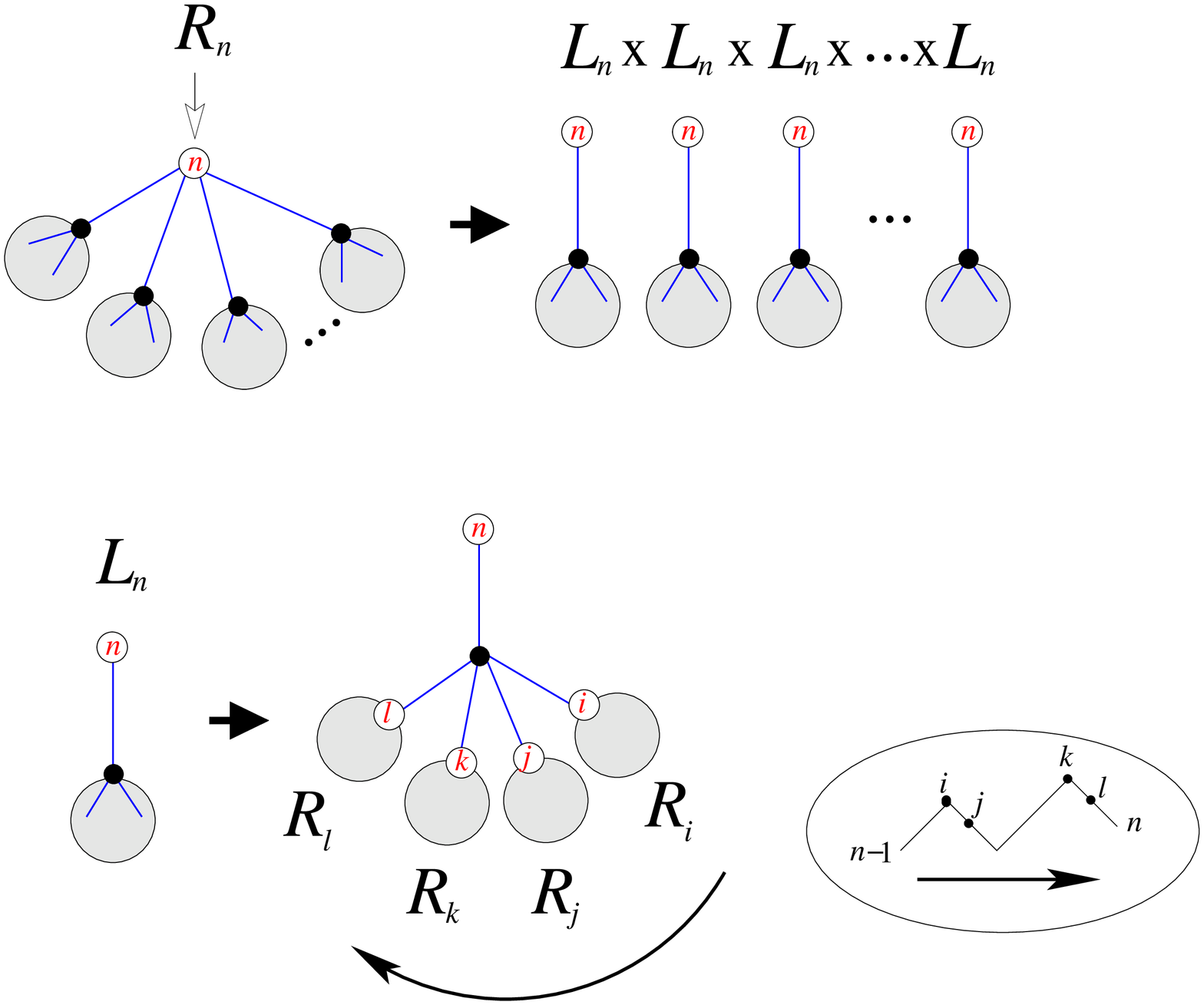}{11.cm}
\figlabel\recursion
Let $R_n\equiv R_n(\{g_{2k}\})$ denote the generating function for rooted
mobiles with root corner labeled $n$ (see figure \recursion). 
Splitting a mobile at its
root vertex leads to an arbitrary number of sub-mobiles, which implies
the relation
\eqn\RtoL{R_n={1\over 1-L_n}\ ,}
where $L_n$ denotes the generating function for mobiles rooted at
a {\it univalent} vertex labeled $n$. 
These new objects may be decomposed according to the sequence of
labels around the unlabeled vertex adjacent to the root. 
By definition of mobiles, such a sequence of length $k$ is in 
correspondence with a walk of length $2k-1$ on the integer line 
with steps of $\pm 1$, starting at position $n-1$ and ending at
position $n$. To each descending step $m\to m-1$ is associated
a mobile rooted at a corner $m$ (see figure \recursion). This implies that
\eqn\recumob{ L_n = \sum_{k=1}^\infty g_{2k} \langle n\vert Q^{2k-1}\vert n-1
\rangle\ ,}
where $Q$ is the operator acting on a formal orthonormal basis 
$|i\rangle$, $i\in \IZ$ (with $\langle j | i \rangle = \delta_{ij}$) as
\eqn\qact{Q |i\rangle=|i+1\rangle +R_i|i-1\rangle\ .}
The operator $Q$ may be understood as a transfer matrix acting {\it
clockwise} around unlabeled vertices.

For instance, the first few valences contribute as
\eqn\somepks{\eqalign{
\langle n|Q|n-1\rangle 
&=1\cr 
\langle n|Q^3|n-1\rangle 
&=R_n+R_{n+1}+R_{n-1} \cr
\langle n|Q^5|n-1\rangle 
& =R_{n+1}R_{n+2}+R_{n+1}R_{n-1}
+R_{n-1}R_{n-2}\cr &\ \ +R^2_{n+1}+R^2_{n} +R^2_{n-1}
+2 R_n(R_{n+1}+R_{n-1})\ .\cr}}
In the case of {\it well-labeled} mobiles, the positivity of labels
is enforced by imposing $R_{-i}=0$ for $i=0,1,2,\cdots$, and by 
considering Eqs. \RtoL-\qact\ for $n\geq 1$ only. With these appropriate
boundary conditions, these equations form a closed set determining
formally the $R_n$'s and $L_n$'s as power series in the $g_{2k}$'s.
The presence of at least one label $1$ is guaranteed by considering
the quantity $R_n-R_{n-1}$, which is therefore identified 
with the generating function for well-labeled mobiles with a distinguished
corner with label $n$. Thanks to our bijection, it is also the generating 
function for bipartite planar maps with an origin and a distinguished 
edge of type $n-1\to n$. As $n\to \infty$, the limit $R$ of $R_n$
corresponds to the generating function for bipartite planar maps with
an origin and a distinguished edge, and satisfies 
\eqn\Requa{R={1\over 1-\sum\limits_{k=1}^{\infty}
{2k-1\choose k} g_{2k} R^{k-1}}\ .}
Alternatively, $R$ is also the generating function for 
(non-necessarily well-labeled) mobiles rooted at a corner labeled,
say, $0$, which has a clear bijective justification.
Similarly, using $\sum_{k\geq 1}(1/k) L_n^k={\rm Log}\, R_n$, we find that
the function ${\rm Log}(R_n/R_{n-1})$ for $n\geq 2$ (resp.\ 
${\rm Log}\, R_1$ for $n=1$) generates 
bipartite planar
maps with an origin and a distinguished vertex at geodesic distance $n$, weighted
by their inverse symmetry factor.

As a final remark, note that equivalent equations have been derived
by the use of an alternative bijection using so-called blossom trees in \
Ref.\ \GEOD, where explicit solutions have been found, due to a remarkable
integrability property. 

\newsec{General case: Eulerian maps with prescribed face degree distribution}

In this section, we extend the construction of section 2 to a 
more general class of maps, namely Eulerian, i.e. face-bicolored maps
with fixed numbers of faces with prescribed color and degree. 
This generalizes the previous case which may indeed be recovered
by imposing that all the faces of a given color, black say, have
degree two, and that all the faces of the other color, white say, have
an even degree, and by contracting all the black faces into single edges.

\subsec{From Eulerian maps to generalized mobiles}

In this section, we show how to associate to an Eulerian map $\cal M$
a labeled tree $\cal T$ that we call a generalized mobile, characterized
by specific label constraints.

\noindent{\underbar{Geodesic distance in an Eulerian map $\cal M$:}}

We start with an Eulerian map $\cal M$ and as before select
an origin vertex. 
Every edge receives a natural orientation
by requiring that the incident black face sits on its right.
We now label each vertex of $\cal M$ by its {\it oriented geodesic 
distance} from the origin, namely the length of any shortest path 
respecting the edge orientations from the origin to this vertex.
Such an oriented path always exists because, starting from an unoriented
path, we may bypass any edge that points in the wrong direction by  
circumventing an incident face. The origin is the only vertex labeled
$0$ and all the labels are non-negative. Moreover, by the
geodesic requirement, the labels $m$ and $n$ on, respectively,
the starting-point and the end-point of any oriented edge
must satisfy $n\leq m+1$.

\noindent{\underbar{Construction rules for $\cal T$:}}
\fig{Rules to construct a tree $\cal T$ out of an Eulerian map $\cal M$.
We first add an unlabeled vertex at the center of each face of the map, with 
the same color (white or black), represented by thick (empty or filled) circles. 
Each oriented edge $e$ of $\cal M$ must connect vertices with labels of the form 
(I) $m\to m+1$ or (II) $m\to n$, with $n\leq m$. 
Accordingly, we draw an edge of $\cal T$ respectively 
connecting (I) the white unlabeled vertex on the left of $e$ to the vertex labeled $m+1$ and
(II) the white and black unlabeled vertices on both sides of $e$. In the latter case
we decorate the edge with labeled flags carrying the labels of the
closest vertices.}{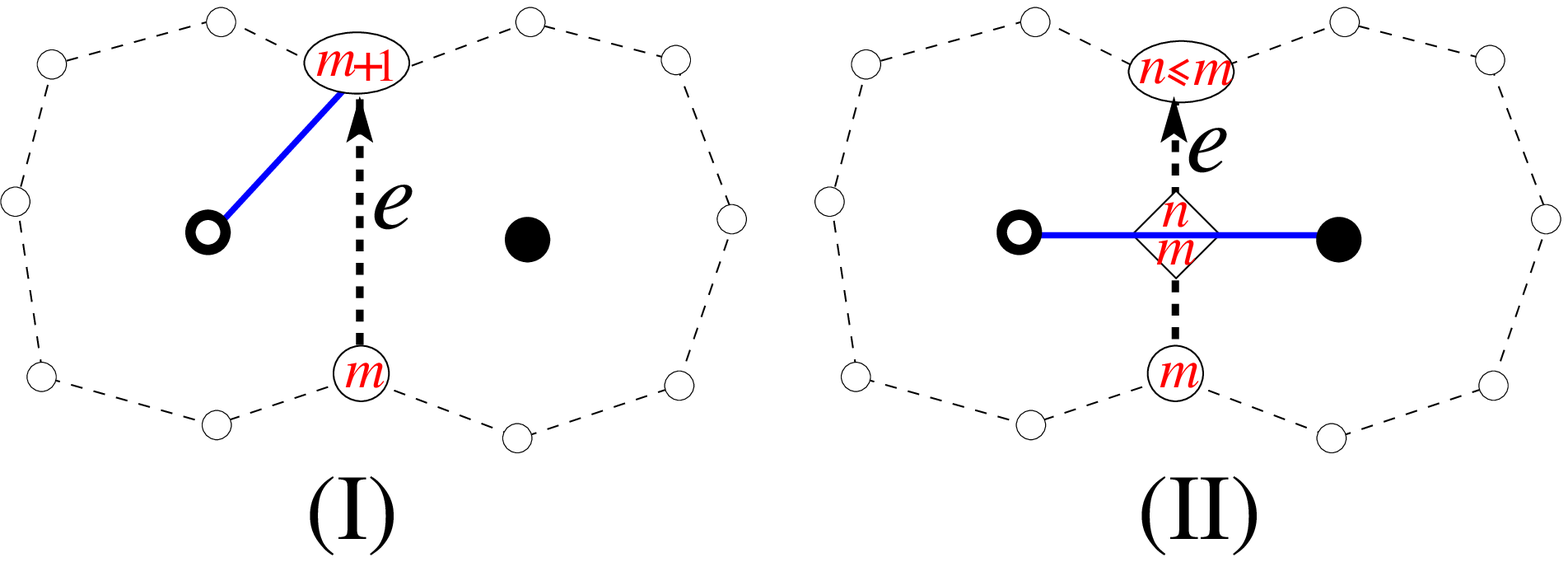}{12.cm}
\figlabel\const
The construction begins by adding a new unlabeled vertex
at the center of each face of $\cal M$. Each unlabeled vertex
receives the color of the face it lies in.
We then consider independently each oriented edge $e$ of $\cal M$ and as
above denote by $m$ and $n$ the labels of its starting- and end-points.
We apply one of the following constructions:
\item{(I)} if $n=m+1$, we draw a new edge from the end-point of $e$ to
the white unlabeled vertex in the white face at the left of $e$
(see figure \const-(I)).
\item{(II)} if $n\leq m$, we drawn a new edge cutting $e$ and linking 
the unlabeled (black and white) vertices in both faces adjacent to $e$.
On each side of this new edge, we add a {\it flag} carrying the label
$m$ or $n$ of the extremity of $e$ lying on this side (see figure \const-(II)). 
\par
After completing this construction for each edge of $\cal M$,
we remove all these original edges. We also erase the origin vertex
because by construction it becomes isolated. We are left with a graph
$\cal T$ made of three types of vertices:
\item{$\bullet$} labeled vertices consisting of all the vertices of $\cal M$
except the origin,
\item{$\bullet$} unlabeled black vertices in correspondence with the black faces
of $\cal M$,
\item{$\bullet$} unlabeled white vertices in correspondence with the white faces
of $\cal M$.
\par 
The edges of $\cal T$ are of two types:
\item{$\bullet$} flagged edges with two flags, one on each side, connecting 
black to white unlabeled vertices
\item{$\bullet$} unflagged edges connecting white unlabeled vertices to 
labeled ones.
\fig{The construction of a tree $\cal T$ out of a sample Eulerian map $\cal M$
with an origin. (a) Each vertex of $\cal M$ is labeled by its oriented geodesic distance from the
origin $0$. (b) We perform the construction of figure \const\ (I)-(II) for each
oriented edge of the map. (c) Erasing the origin and the original edges of $\cal M$ produces
a tree $\cal T$ carrying labels on some vertices and edges.}{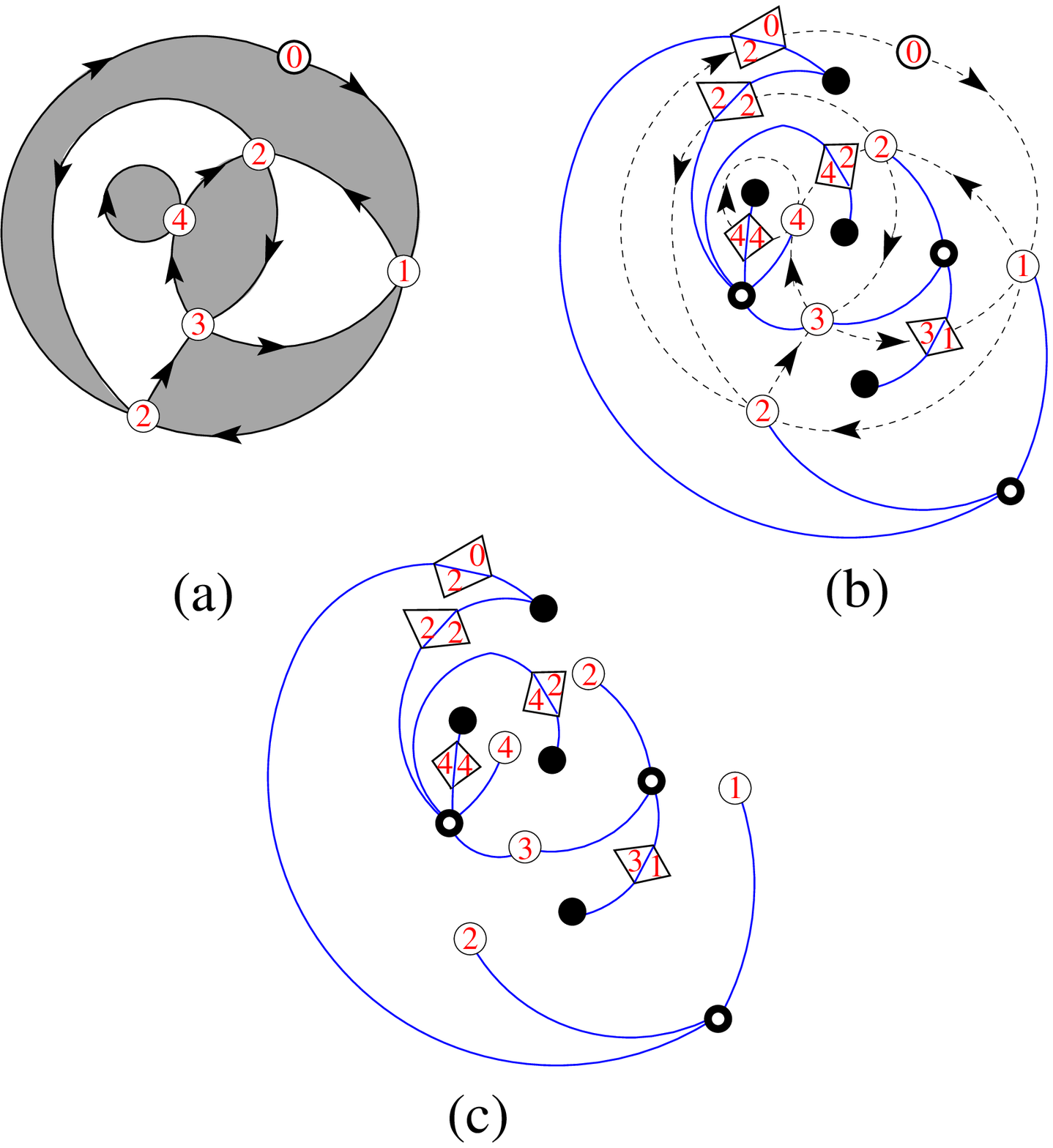}{11.cm}
\figlabel\eulermaptomob
Figure \eulermaptomob\ shows an example of our construction on a
sample Eulerian map.
\par 
\noindent{\underbar{Proof that $\cal T$ is a tree:}}

As before, $\cal T$ is a plane tree: it has no cycles and a single
connected component. For the first statement, assume by contradiction
the existence of a cycle, delimiting an interior domain not containing
the former origin. We now consider the smallest label $n$ among
\item{-} the labels of the labeled vertices belonging to the cycle,
\item{-} the labels on the flags along the cycle and lying in the interior.
\par
If $n$ is attained for a labeled vertex on the cycle, we deduce
from the construction (I) that there is a vertex labeled $n-1$ in the interior.
Consider an oriented geodesic path on $\cal M$ to
this latter vertex: it necessarily intersects the cycle at a labeled vertex 
since all its edges are of type (I) above hence they never cut an
edge of $\cal T$. The label at the intersection is strictly smaller than $n$.
Otherwise, if $n$ is not attained for a labeled vertex, but only by flags,
we deduce from the construction (II) the existence of vertex labeled $n$ in
the interior, hence by considering again a geodesic path to this 
latter vertex, the existence of a vertex with label $\leq n$ on the cycle.
In both cases, we have a contradiction.

The connectedness is proved as before by a counting argument:
since $\cal T$ contains no cycle, its number $c$ of connected
components is the difference between its number of vertices and
its number of edges. Denoting by $V$, $E$, and $F$ the numbers
of vertices, edges and faces of the original Eulerian map $\cal M$,
the number of edges of $\cal T$ is obviously $E$ by construction,
while the number of its vertices is $V+F-1$; hence $c=1$ by the Euler 
relation for $\cal M$. This concludes the proof that $\cal T$
is a plane tree.

\noindent{\underbar{Label characterization of $\cal T$:}}
\fig{Constraints $P_{\bullet}$
on labels around a black unlabeled vertex. Clockwise around the vertex
two successive (flag) labels $m$ and $n$ must either satisfy $n\leq m$
if the flags share the same edge, or $n\geq m$ otherwise. In the latter case,  
$p=n-m+1$ may be interpreted as the number of labeled vertices incident in $\cal M$
to the corresponding black face, and in between the two flags.}{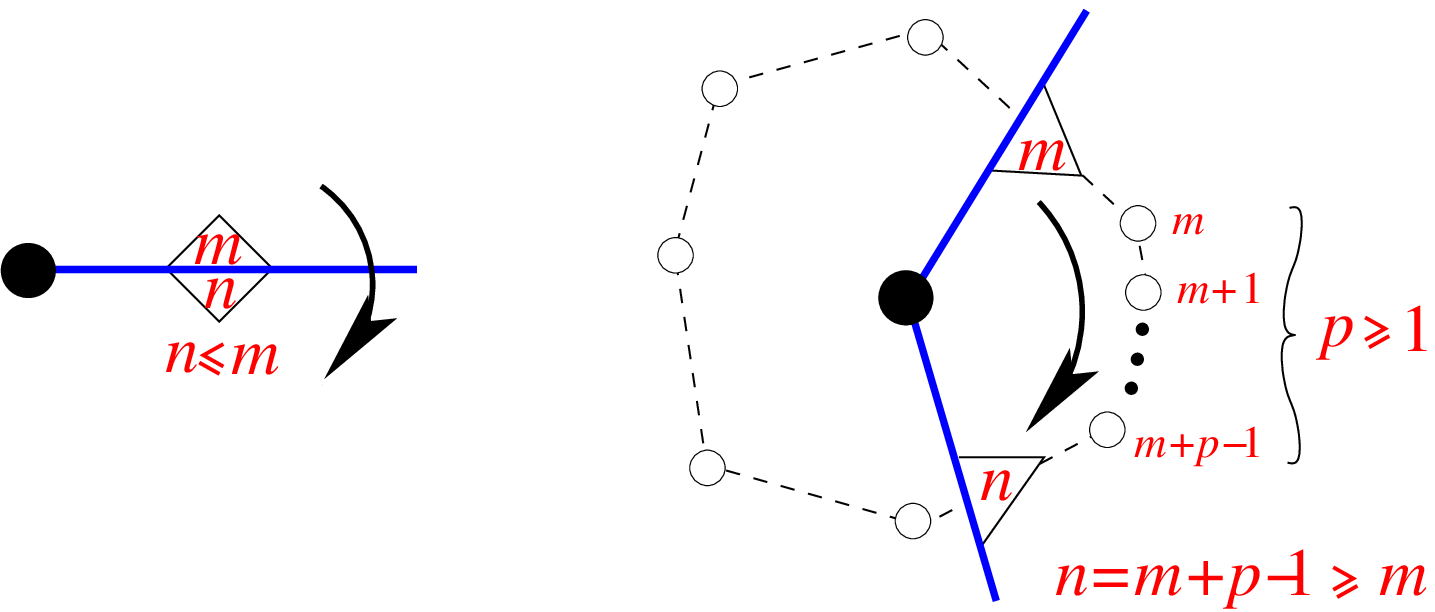}{10.cm}
\figlabel\eulerconstblack
We now investigate the properties of labels in $\cal T$
inherited from those of $\cal M$.
Each black unlabeled vertex is incident to flagged edges
only. 
The flag labels have the following property:
\item{(P$_\bullet$)}
for each black unlabeled vertex of $\cal T$, the sequence of flag labels
of its incident edges, when read clockwise, is non-increasing at 
each edge crossing and non-decreasing between two consecutive edges 
(see figure \eulerconstblack).
\par

\fig{Constraints $P_\circ$ on labels around a white unlabeled vertex read clockwise.
Two successive flag labels $m$ and $n$ sharing the same edge must satisfy $n\geq m$.
For consecutive edges,
a vertex labeled $m$ is followed by a label $m-1$ (either on a vertex or a flag),
while a flag labeled $m$ is followed by a label $m$ (either on a vertex 
or a flag).}{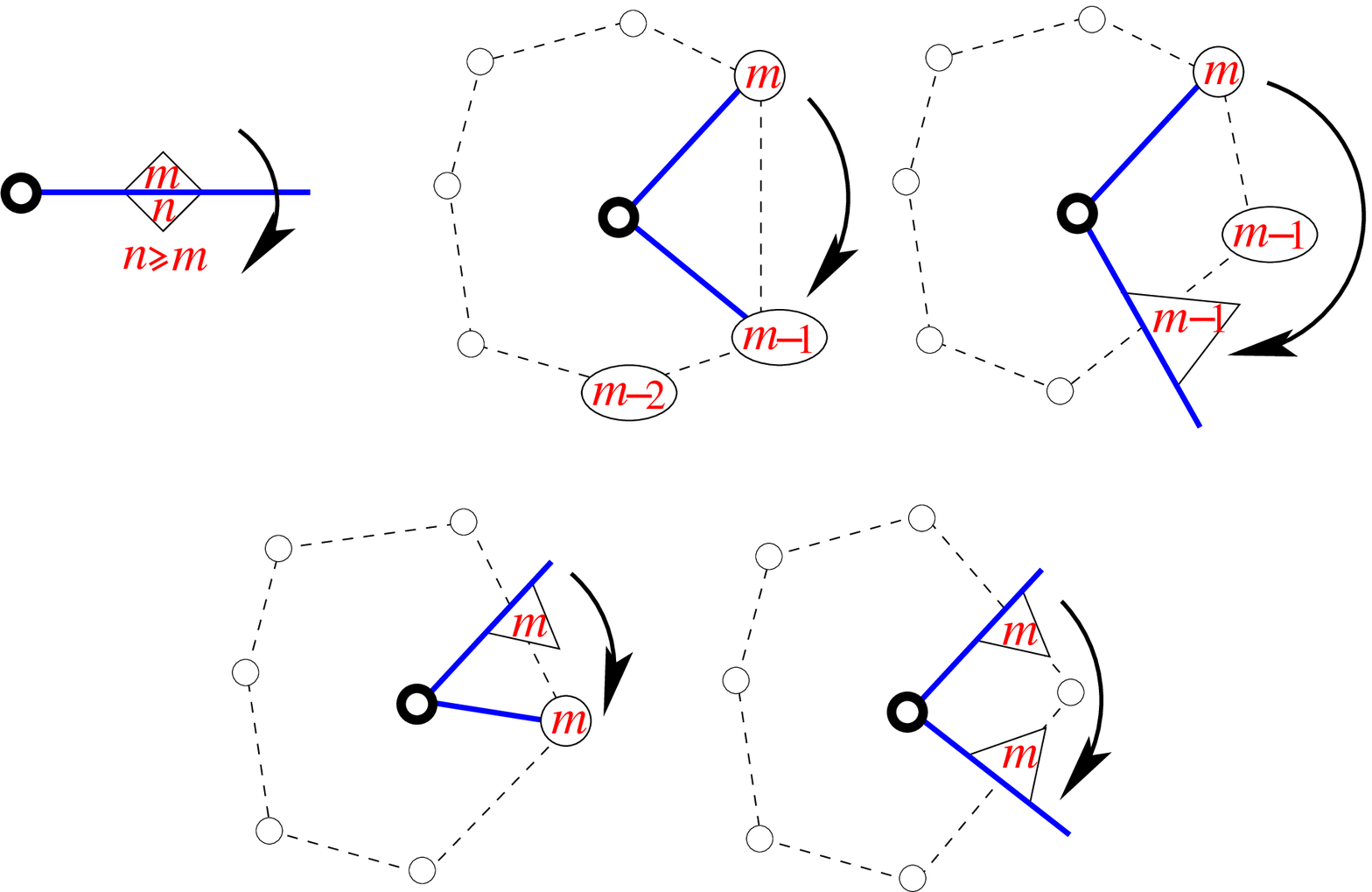}{13.cm}
\figlabel\eulerconst

\noindent Each white unlabeled vertex is incident to both flagged edges
and unflagged edges leading to labeled vertices. The labels
have the following property:
\item{(P$_\circ$)}
for each white unlabeled vertex, the sequence of labels on incident
flagged edges and on adjacent labeled vertices, when read clockwise,
is non-decreasing at each edge crossing, decreasing by $1$ after
each labeled vertex and stationary between a flagged edge and the
next label (see figure \eulerconst). 
\par
\fig{ The coding of the sequence of labels around white (left)
and black (right) faces of $\cal M$ (top) and unlabeled vertices
of $\cal T$ (bottom) via a closed walk (center) with 
increments $\geq -1$.  The successive
labels around faces of $\cal M$ when read clockwise match the successive
heights of the walk when read from left to right (white faces) or
right to left (black faces).
The vertex labels around an unlabeled white vertex of $\cal T$   
match the heights preceding steps of $-1$ on the walk, while flag
labels match all other steps, of non-negative integer height.
The flag labels around an unlabeled black vertex of $\cal T$ 
match the heights of all steps of non-positive 
integer height.}{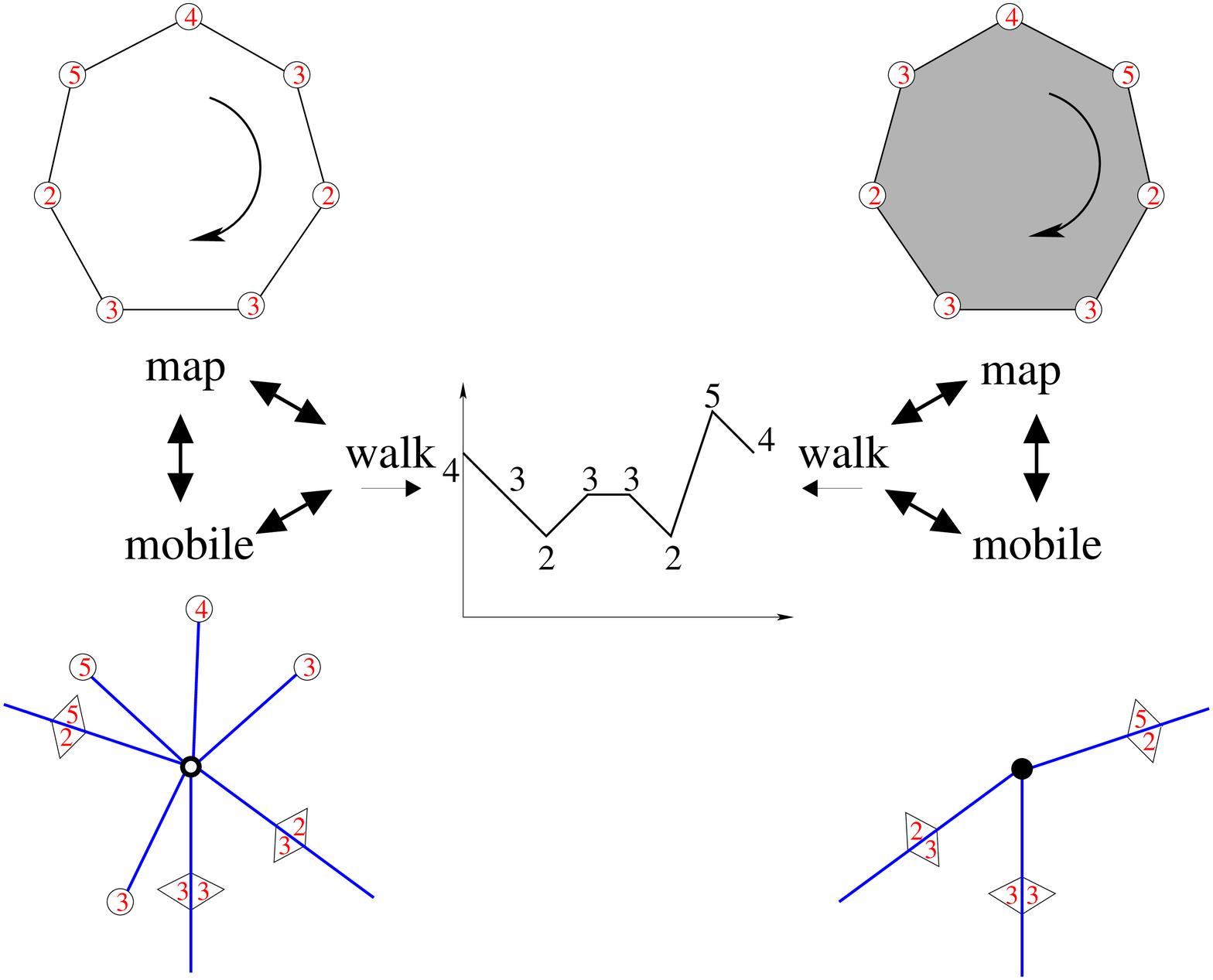}{13.cm}
\figlabel\walks
\noindent $P_\bullet$ and $P_\circ$ simply rephrase the 
properties of labels around faces of $\cal M$, namely that the 
clockwise sequence of labels around a black (resp.\ white)
face are in one-to-one correspondence with closed walks with ascending 
steps  of arbitrary (possibly zero) integer height, and descending steps of $-1$. 
The walk must be read from left to right (resp.\ right to left) 
in the case of a white (resp.\ black) face, as indicated in figure \walks.
Note also that the valence of the faces of $\cal M$ is nothing but
the length (number of steps) of the corresponding walks.
As shown in figure \walks, the coding of a walk in terms of tree labels
is slightly different around black and white vertices.

We define a generalized mobile to be a plane tree such that:
\item{(i)} its vertices are of three types: black unlabeled vertices,
white unlabeled vertices and labeled vertices carrying integer labels,
\item{(ii)} its edges are of two types: flagged edges with one flag
on each side carrying an integer label, which connect unlabeled 
vertices of different colors, and unflagged edges, each connecting a labeled
vertex to a white unlabeled one,
\item{(iii)} the labels obey the properties (P${}_\bullet$) and
(P${}_\circ$) above.
\par
\noindent If the mobile moreover satisfies the additional rule:
\item{(iv)} all vertex labels are strictly positive, all 
flag labels are non-negative, and there is at least a flag labeled
$0$,
\par
\noindent then it will be called {\it well-labeled}.

The above construction provides a mapping from Eulerian planar maps
to well-labeled generalized mobiles. This mapping is a bijection
whose inverse is constructed in the next section.

\subsec{Inverse construction}

We start from a well-labeled generalized mobile $\cal T$. As before, we define
a corner of $\cal T$ as a sector with apex at a labeled vertex and
delimited by two consecutive (unflagged) incident edges. Each corner
receives the label of its apex. We now define successors for both 
corners and flags:
\item{-} for each corner $C$ with label $n\geq 2$, its successor $s(C)$ 
is the first encountered corner with label $n-1$ when going clockwise
around the tree,
\item{-} for each flag $F$ with label $n\geq 1$, its successor $s(F)$ is 
the first encountered corner with label $n$ when going clockwise
around the tree.
\par
\noindent The existence of successors is a consequence of the
rules (iii) and (iv) defining well-labeled generalized mobiles. 
Indeed, the sequence of labels read clockwise around $\cal T$ 
is non-decreasing after a flag, and decreases by one after a corner. 
Hence between a corner with label $n\geq 2$ and a flag labeled $0$, 
there is at least one label $n-1$ after which the sequence decreases, thus
corresponding to a corner. Similarly between a flag with label $n\geq 1$ 
and a flag labeled $0$, there is at least one label $n$ after which the 
sequence decreases, again corresponding to a corner.

We construct the map $\cal M$ associated to $\cal T$ by connecting
each corner with label $\geq 2$ and each flag with label $\geq 1$
to its successor, which can be done in such a way that these newly
created edges do not intersect. This is possible because of the 
nested structure of successors, a consequence of the above property
of the label sequence  around $\cal T$. At this point,
all the corners labeled $1$ and flags labeled $0$ are in the same external
face and we connect them all to a new origin vertex labeled $0$ inside
this face.
Finally we erase all the edges of $\cal T$, whether flagged or unflagged,
as well as all the unlabeled vertices, whether black or white. 
The result is a map $\cal M$ with an origin which is
connected as each vertex is connected to the origin via a chain
of successors.
\fig{The construction of an Eulerian map $\cal M$ with an origin out of a sample 
well-labeled generalized mobile $\cal T$.
Starting from the mobile (a) we connect each labeled corner and flag to its successor
via an oriented edge
(dashed arrow). This includes adding a vertex labeled $0$, successor of
all corners labeled $1$ and flags labeled $0$.  Erasing all unlabeled
vertices and all edges of the mobile produces an Eulerian planar map with an origin
(c). Its vertex labels moreover encode their geodesic distance 
from the origin.}{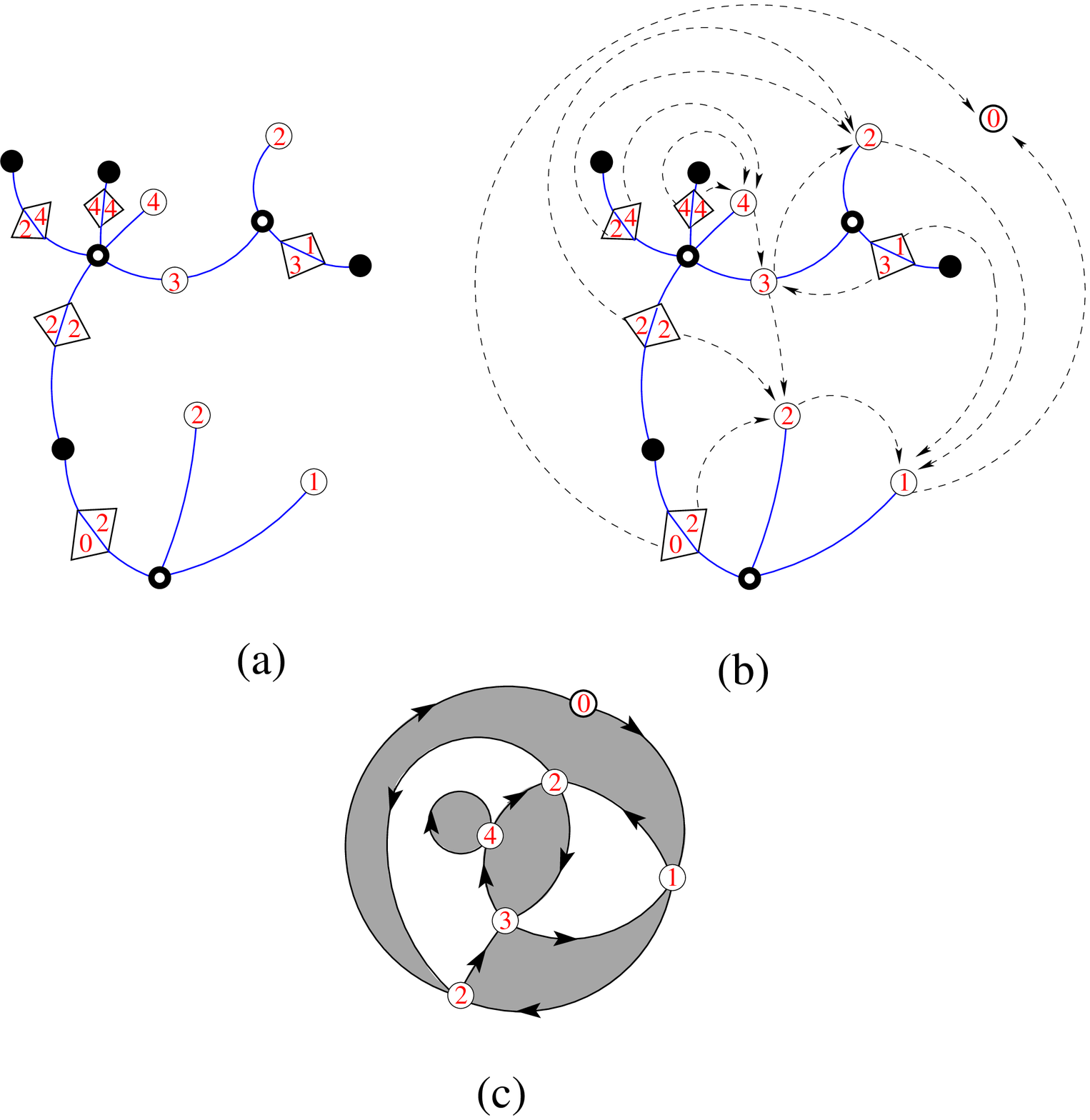}{11.cm}
\figlabel\mobtoeulermap

With the help of figure \mobtoeulermap, the reader can easily check 
that $\cal M$ is indeed an Eulerian map,
and that its vertices are labeled by their oriented geodesic distance from
the origin. Moreover, the above construction and that of Sect.\ 3.1
are the inverse of one another. This can be proved by following the 
same line of arguments of in Sect.\ 2.3, suitably generalized so as to
include the flagged edges. We omit the details here.


\subsec{Generating functions}

As in section 2.4, we may derive recursion relations for 
the generating functions for mobiles, also interpreted
as generating functions for Eulerian maps. Again we
attach weights $g_k$ (resp.\ ${\tilde g}_k$) to white
(resp.\ black) $k$-valent faces of the Eulerian map.
\fig{The generating functions $R_n,L_n,B_{m,n},W_{m,n}$ for respectively
rooted mobiles with root corner labeled $n$, rooted mobiles with {\it univalent} root
vertex labeled $n$, and half-mobiles with unlabeled univalent root vertex
attached to a flagged edge with labels $n,m$ incident to a black, white vertex. }{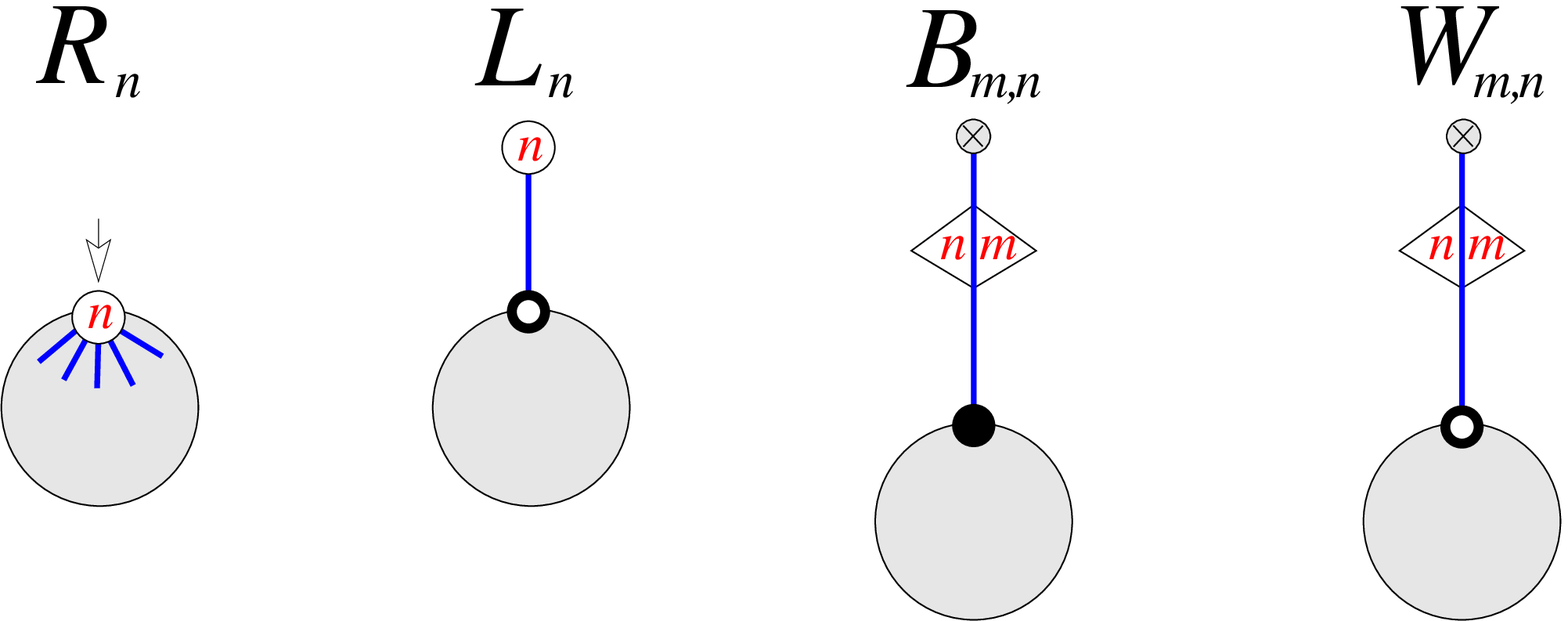}{10.cm}
\figlabel\subtrees
More precisely, let $R_n$ denote the generating function
for rooted mobiles, i.e. with a distinguished corner labeled $n$.
As in section 2.4, by splitting such a mobile at its root vertex, 
we obtain the equation
\eqn\RtoLagain{R_n={1\over 1-L_n}\ ,}
where $L_n$ denotes the generating function for mobiles
rooted at a univalent vertex labeled $n$.
To proceed further in the decomposition of these objects,
we are led to introduce new classes of objects, namely {\it half-mobiles},
obtained by cutting a mobile across a flagged edge.
These are represented in figure \subtrees\ as trees starting with a
univalent unlabeled root vertex incident to a flagged edge.
The rules $P_\circ$ and $P_\bullet$ are obeyed around each
unlabeled vertex except the root. We denote by $W_{m,n}$ (resp.\  
$B_{m,n}$)
the generating function for half-mobiles whose root vertex is adjacent 
to a white (resp.\ black) unlabeled vertex via a flagged edge with flag 
$m$ on the right and $n$ on the left when going towards the root
(see figure \subtrees). By the rules $P_\circ$ and $P_\bullet$,  
$W_{m,n}=0$ for $m<n$, while $B_{m,n}=0$ for $m>n$.

By decomposing a mobile rooted at a univalent vertex labeled $n$
around its adjacent unlabeled white vertex, we obtain the recursive relation
\eqn\LtoB{L_n=\sum_k g_k \langle n | Q^{k-1} | n-1 \rangle\ ,}
where $Q$ is the operator acting on a formal orthonormal basis as
\eqn\Qacttwo{Q|i\rangle = \sum_{j\geq i} B_{i,j}|j\rangle + R_i|i-1\rangle\ ,}
which generates the sum over all walks from height $n-1$ to height $n$
with steps $\geq -1$ and with weights
\item{-} $R_i$ per step $i \to i-1$
\item{-} $B_{i,j}$ per step $i \to j$ with $j \geq i$
\item{-} $g_k$ per walk of length $k-1$ (corresponding to a $k$-valent
white face).
\par
\noindent This follows from the correspondence displayed in figure
\walks\ between white vertex neighborhoods in mobiles, and walks : the
operator $Q$ may be understood as a transfer matrix acting
clockwise around white vertices and therefore corresponds to reading
walks from left to right in figure \walks.

Similar recursive relations for $W_{m,n}$ with $m\geq n$ follow from
the same considerations, leading to
\eqn\WtoB{W_{m,n}=\sum_k g_k \langle n | Q^{k-1} | m \rangle\ .}

Finally, repeating the exercise around black vertices leads to the
following recursive relation valid for $m\leq n$
\eqn\BtoW{B_{m,n}=\sum_k \tilde{g}_k \langle n | \tilde{Q}^{k-1} | m \rangle\ ,}
where $\tilde{Q}$ acts as
\eqn\tQact{\tilde{Q}|i\rangle = |i+1\rangle + \sum_{j\leq i} W_{i,j}|j\rangle\ ,}
and may be interpreted as a transfer matrix generating black vertex
neighborhoods in the clockwise direction, or, equivalently, walks (read from
right to left in figure \walks) from height $m$ to $n$, with steps
$\leq 1$ and weights
\item{-} $1$ per step $i \to i+1$ (no connected subtree)
\item{-} $W_{i,j}$ per step $i \to j$ with $j \leq i$
\item{-} $\tilde{g}_k$ per walk of length $k-1$ (corresponding to a
$k$-valent black face).
\par
\noindent 
For well-labeled mobiles these equations are only valid for strictly
positive vertex labels and non-negative flag labels with the boundary
conditions
\eqn\BCeul{R_i=L_i=0, \quad {\rm for}\ i \leq 0\ ,
\qquad W_{i,j}=B_{i,j}=0, \quad {\rm for}\ i<0\ {\rm or}\ j<0\ .}
With these boundary conditions, the set of relations \RtoLagain, \LtoB, 
\WtoB\ and \BtoW\
form a closed set determining order by order in the $g$'s
and $\tilde{g}$'s all the generating functions involved. 

As before, the function $R_n-R_{n-1}$ generates well-labeled mobiles
with root corner labeled $n$, and with at least one label $1$. Mapwise, 
this is nothing but the generating function for Eulerian planar maps with
an origin vertex and a distinguished oriented edge $n-1\to n$. Moreover,
for $m\geq n$, the function $B_{n,m}W_{m,n}-B_{n-1,m-1}W_{m-1,n-1}$
generates the Eulerian planar maps with an origin vertex and a distinguished
oriented edge $m\to n$.
In the limit of large labels, we have $R_i \to R$, $L_i \to L$,
$W_{i,i-m}\to W_{m}$ and $B_{i,i+m}\to B_{-m}$ ($m\geq 0$) and  
we recover the relations derived in Ref. \BMS\ (up to 
some rescalings). The function $B_{-m}W_{m}$ is the
generating function for Eulerian planar maps with an origin vertex
and a distinguished edge $k\to k-m$ for some $k$.

As before, the function ${\rm Log}(R_n/R_{n-1})$ for $n\geq 2$ 
(resp.\ ${\rm Log}\, R_1$ for $n=1$) generates Eulerian planar maps 
with an origin and a distinguished vertex at geodesic distance $n$, weighted by their
inverse symmetry factor.

As a final remark, we note that the operators $Q$ and $\tilde{Q}$ 
of eqns \Qacttwo\
and \tQact\ are identical to those used in the orthogonal polynomial
solution of the two-matrix model generating Eulerian maps,
respectively implementing the effect of multiplication by an
eigenvalue of either matrix. Moreover the recursion relations above
are very similar to those governing the genus expansion of this matrix
model with the only difference that eq \RtoLagain\ should be replaced
by $R_n=(n/N)/(1-L_n)$ for matrices of size $N \times N$.

\newsec{Some particular cases of interest}

In this section we show how the construction of Sect.\ 3 simplifies for
two particular classes of planar maps, namely the so-called
$p$-constellations \CONST\ and the general planar maps with arbitrary
(even or odd) valences \CENSUS.

\subsec{$p$-constellations}

A $p$-constellation is an Eulerian planar map
such that all black faces have degree $p$
while white faces may have degrees multiples of $p$.
Note that the case of maps of even degrees of Sect.\ 2
corresponds to constellations with $p=2$ by contracting all $2$-valent 
black faces into single edges.

As before, the vertices of a $p$-constellation with a fixed origin vertex
are labeled as in an ordinary Eulerian map by their oriented geodesic distance 
from the origin. However these labels are more constrained in a 
$p$-constellation: {\it the labels along each oriented edge may
either increase by $1$ or decrease by $p-1$}.
This is proved by noticing that the length of any two oriented paths from the
origin to a given vertex are congruent modulo $p$. 
Indeed, the difference between the two lengths is equal to the sum of
the perimeters of the black faces in between the two paths minus
that of the white ones, all of which are multiples of $p$.
As a consequence, the labels along oriented edges must increase 
by $1$ modulo $p$, and can neither increase by more than 
$1$ from the geodesic requirement nor decrease by more than $p-1$ 
because we may use the bypass around the incident $p$-valent black face.

This property results in a drastic simplification of the $p$-mobiles 
built out of $p$-constellations from the construction of Sect.\ 3.1. 
Indeed each black face of a $p$-constellation becomes a 
{\it univalent} black unlabeled vertex in the associated $p$-mobile. 
This is readily seen by noting that the clockwise sequence of the 
$p$ labels around a black face is a walk with $(p-1)$ (ascending) steps
of $+1$ and {\it one} (descending) step of $-(p-1)$. 
This unique descending step corresponds to the unique edge
incident to the corresponding unlabeled black vertex.
This edge carries flags with labels differing by $(p-1)$. This exhausts
all flagged edges of the $p$-mobile. 

Turning to generating functions, the $p$-constellations are enumerated with
a weight $g_{pk}$ per $pk$-valent white face ($k\geq 1$) and 
a weight ${\tilde g}_p$ per black face. As in Sect.\ 3.3, the enumeration 
is performed via the bijection to well-labeled $p$-mobiles. 
Introducing as before the generating functions $R_n$ and $L_n$ for rooted $p$-mobiles,
and those
for half-mobiles $B_{m,n}$ and $W_{m,n}$,
the above restrictions imply that
\eqn\bmnrestrict{B_{m,n}={\tilde g}_p \delta_{n,m+p-1}} 
for the trivial
half-mobiles whose root vertex is attached to a black univalent vertex. 
The recursive relations \RtoLagain\ and \LtoB\ remain the same,
but with \Qacttwo\ replaced by
\eqn\qandq{Q|i\rangle = {\tilde g}_p |i+p-1\rangle+ R_{i} |i-1\rangle\ .}
Equations \RtoL, \recumob\ and \qandq\ form a closed set of equations 
determining $R_n$ order by order in the $g$'s 
and ${\tilde g}$'s, while the $W_{m,n}$ are determined by \WtoB.

Note that the black vertices can be removed without loss of information,
leading to a slightly simpler definition of $p$-mobiles, now made of only white unlabeled and labeled
vertices, with edges connecting only labeled to unlabeled ones, and with the
constraint that the clockwise labels around any unlabeled vertex must either decrease by $1$
or increase by a quantity of the form $k(p-1)-1$, $k=1,2,3,...$

The above discussion includes the particular case $p=2$ of 
bipartite planar maps, and 
the simplified definition of $2$-mobiles matches 
exactly that of the mobiles of Sect.\ 2.
Moreover, upon taking ${\tilde g}_2=1$, the above equations directly 
reproduce eqns. \RtoL, \recumob\ and \qact.

\subsec{Maps with arbitrary valences}

The case of arbitrary planar maps with prescribed face valences may be seen 
as yet another particular case of Eulerian planar maps, by imposing 
the condition that all the black faces have valence $2$.
These faces may indeed be contracted into regular edges of the ordinary planar map.
As a result, these edges may be traversed in both directions; hence we now deal
with an unoriented geodesic distance.
\fig{A sample planar map with arbitrary valences and an origin. 
The edges are inflated into two-valent black faces producing 
an Eulerian planar map with only two-valent black faces.
As usual, vertices are labeled
by their geodesic distance from the origin. The orientation of the edges 
becomes irrelevant.}{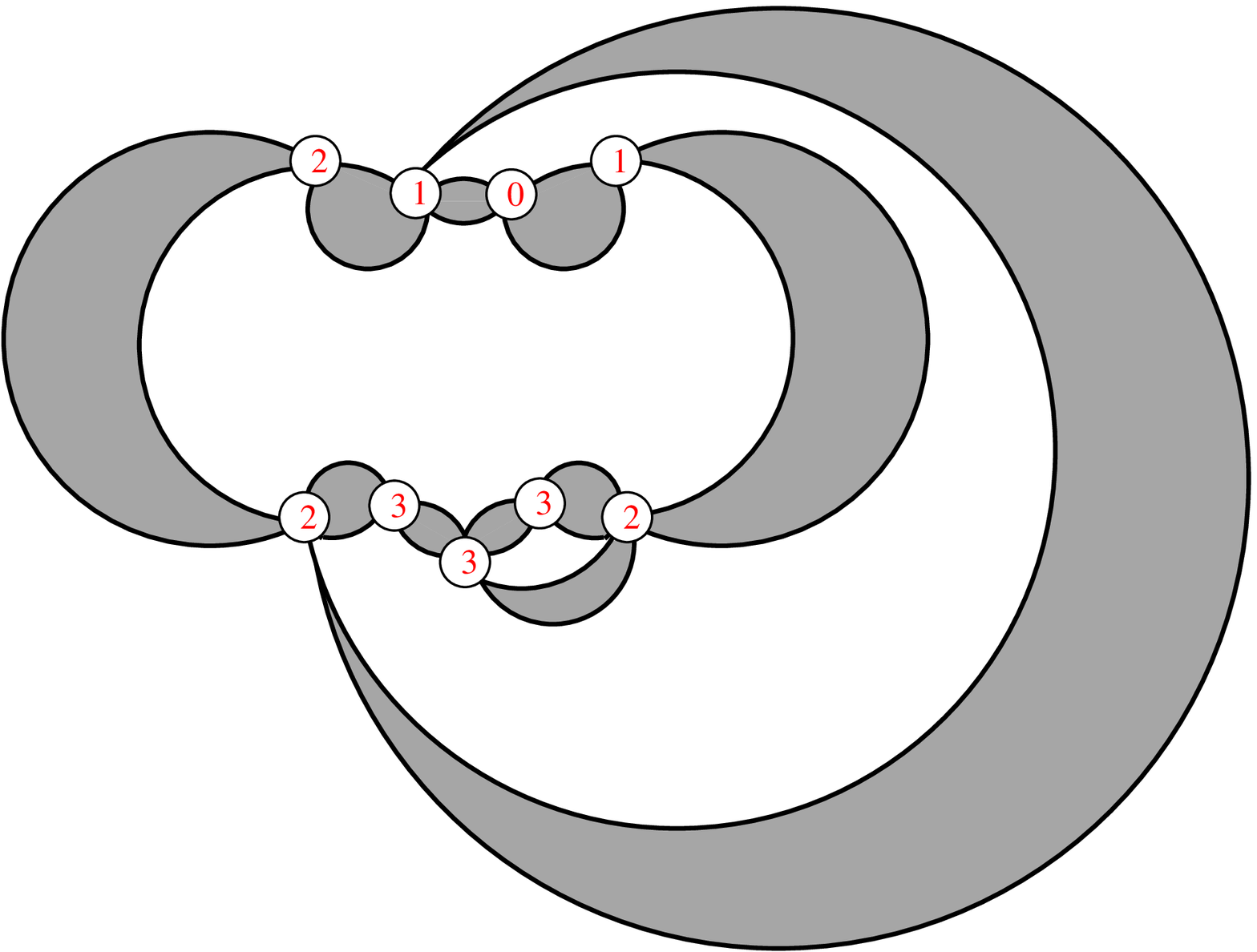}{7.cm}
\figlabel\arbitmap
Again this restricts the label configurations:  
{\it the labels along each oriented edge may
differ by $+1,0$ or $-1$} (see figure \arbitmap). 

\fig{The two possible situations around an unlabeled black vertex
associated with a two-valent black face in an Eulerian map. 
The arrows indicate simplification rules for the mobiles associated to
planar maps with arbitrary valences.}{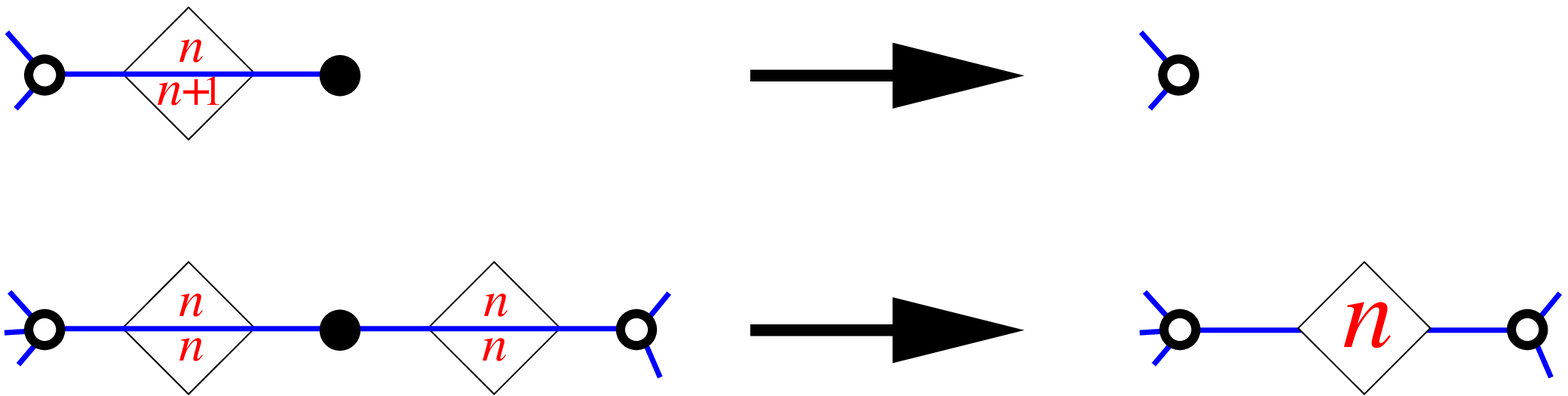}{10.cm}
\figlabel\simplif
On the associated mobile, we have only two types of unlabeled black vertices
as shown in figure \simplif:
\item{$\bullet$} univalent black vertices, with an incident flagged edge,
labeled $n,n+1$ 
\item{$\bullet$} bivalent black vertices, with two incident flagged edges, both
labeled $n,n$.

Turning to the generating functions, we immediately obtain, in addition to eqs. \RtoL\
and \LtoB
\eqn\bmnarbit{ B_{n,n+1}=1\ , \qquad B_{n,n}= W_{n,n}\ ,}
where we simply have restricted ourselves to ${\tilde g}_k=\delta_{k,2}$,
while 
\eqn\wmnarbit{ W_{n,m}=\sum_k g_k \langle n | Q^{k-1} | m \rangle\ ,\qquad m=n,n+1\ ,}
where the operator $Q$ reduces to
\eqn\Qarbit{Q|i\rangle = |i+1\rangle+ W_{i,i}|i\rangle+R_i|i-1\rangle\ .}

\fig{A sample planar map with arbitrary valences and an origin together
with its associated (simplified) mobile.}{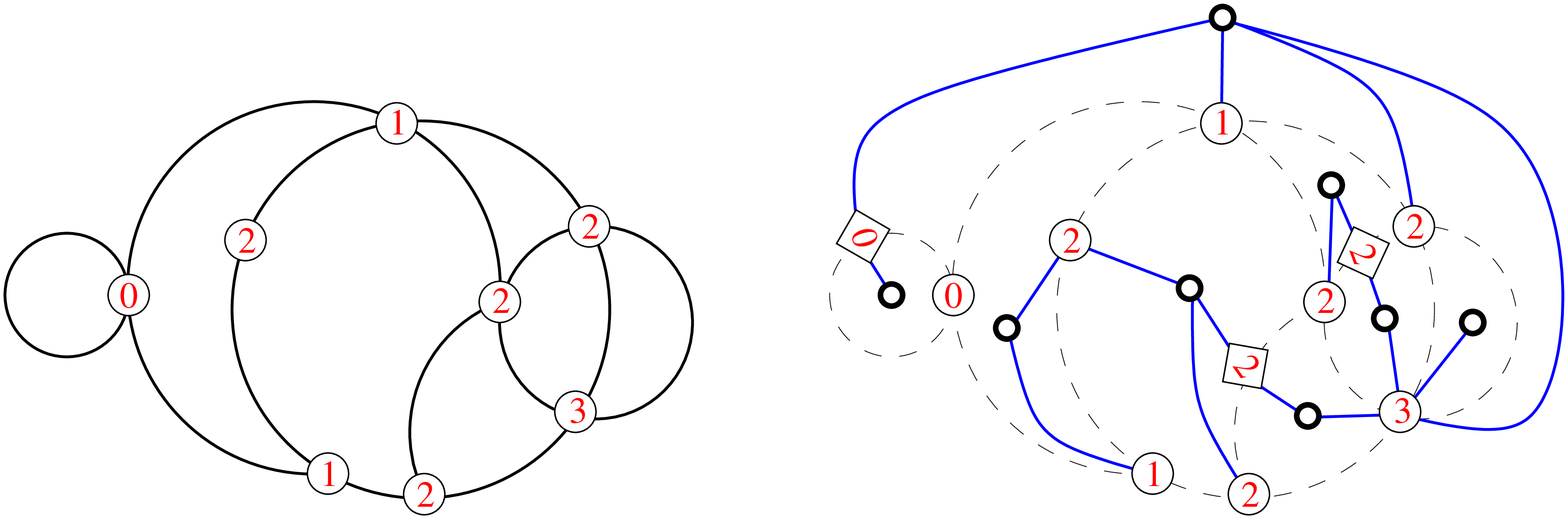}{13.cm}
\figlabel\simpmap

The corresponding mobiles may be simplified without loss of information
by erasing all univalent black vertices and replacing the bivalent ones 
and their two incident flagged edges (with, say, labels $n$) by a new type of edge, 
carrying a unique flag (labeled $n$), as illustrated in figure \simplif.  
The simplified mobiles therefore have only white unlabeled and labeled vertices
with two types of edges: unflagged ones connecting labeled to unlabeled vertices,
and flagged ones connecting unlabeled vertices to one another. The labels obey the following rule
clockwise around unlabeled vertices: a vertex labeled $n$ is followed by a vertex or flag labeled
at least $n-1$, while a flag labeled $n$ is followed by a vertex or flag labeled at least $n$.
Denoting by $S_n=W_{n,n}$ the generating function for half-mobiles with an incident
edge flagged $n$, we get the system of equations 
\eqn\arbiteqs{\eqalign{ R_n&={1\over 1-\sum_{k\geq 1}g_k \langle n|Q^{k-1} |n-1\rangle} \cr
S_n&=\sum_{k\geq 1} g_k \langle n|Q^{k-1} |n\rangle \cr
Q|i\rangle &= |i+1\rangle +S_i |i\rangle +R_i |i-1\rangle\ .\cr}}
Again, this determines $R_n$ and $S_n$ completely order by order in the $g$'s.
As before, the function $R_n-R_{n-1}$ generates arbitrary planar maps with an origin and
a distinguished edge $(n-1,n)$, while $S_n^2-S_{n-1}^2$ generates arbitrary planar maps with an origin and
a distinguished edge $(n,n)$.

\newsec{Conclusion}

In this paper, we have established a general bijection between Eulerian
planar maps with an origin vertex and so-called mobiles, namely
plane trees with specific decorations. The latter carry in particular the
information of the oriented geodesic distance from the origin on the
original map, and allow for an easy enumeration of maps with prescribed
face valences and various distinguished vertices or edges at fixed geodesic 
distance from the
origin. We conclude by discussing some interesting applications.

 \noindent{\underbar{Local environment in infinite maps:}}

The bijection gives access to the statistics of the {\it local environment}
of a vertex (the set of vertices geodesically close to that vertex)
in large Eulerian planar maps, by extracting the singular parts of
the various generating functions involved (see Ref.\ \ONEWALL\ for a general
scheme and examples). For instance we may extract from the singular part
of $R_n-R_{n-1}$ the average number of edges at geodesic distance $n$
(i.e. of the form $(n-1,n)$) from a fixed origin in infinite, say,
quadrangulations, i.e. bipartite planar maps with only tetravalent faces,
and obtain the equation:
 \eqn\quadrangaverage{ \langle e_{n-1,n} \rangle= {6\over 35} {
(n^2+2n-1)(5n^4+20n^3+27n^2+14n+4) \over n(n+1)(n+2) }\ .}
 Similarly, we may also extract from the singular part of ${\rm
Log}(R_n/R_{n-1})$ the average number of vertices at geodesic distance $n$
from an origin in infinite quadrangulations, with the result
 \eqn\verquadinfinite{ \langle v_{n} \rangle= {3\over35}
\big((n+1)(5n^2+10n+2)+\delta_{n,1}\big)\ .}
 \medskip
 \noindent{\underbar{Integrability properties:}}

The recursion relations obtained in this paper match those found in 
Ref.\ \GEOD\ through the bijection with blossom trees, but now with a different
(dual) interpretation in terms of mobiles. Many of these equations were
found to be integrable (i.e. have a maximal number of independent
integrals of motion) and have exact solutions as given in \GEOD, which
directly translate into compact expressions for the generating functions
of rooted mobiles, as well as for the corresponding rooted maps.
 \medskip
 \noindent{\underbar{Branching processes:}}

The mobiles introduced in this paper are (decorated) labeled trees.  It is
tempting to interpret them as spatially branching processes, by viewing
the tree as coding the genealogy of a population while the labels code
positions of its members on the integer line. Such an interpretation was
exploited in Ref.\ \LALLER\ to relate the generating function $R_n$ for
planar quadrangulations to the escape probability of a species from a
domain of the integer line. It would be interesting to see whether the
rules for labels in mobiles have a natural translation in terms of
population evolution, and to study their scaling limit around 
multicritical points.

\bigskip
\noindent{\bf Acknowledgments:} 
All the authors acknowledge the support of the EU network on 
``Discrete Random Geometry", grant HPRN-CT-1999-00161.

\listrefs
\end